\documentclass{article}
\usepackage{amsthm}
\usepackage{xcolor}
\usepackage[english]{babel}
\newcommand{\reals}{{I\kern-.35em R}}
\newcommand{\Reals}{\overline{I\kern-.35em R}}
\newtheorem{theor}{Theorem}[section]

\newtheorem{prop}[theor]{Proposition}
\newtheorem{lema}[theor]{Lemma}
\newtheorem{defi}[theor]{Definition}

\newcommand{\co}{\overline{co}\hbox{ }}

\def\ball{\overline{I\kern -.35em B}}

\def\tto{\;{\lower 1pt \hbox{$\rightarrow$}}\kern -12pt
\hbox{\raise 2.8pt \hbox{$\rightarrow$}}\;}

\usepackage{amsmath,amsxtra, amssymb, latexsym, amsfonts, amscd } 
\usepackage{multicol,color}
\usepackage{bbm}
\usepackage{float}
\usepackage{soul}
\usepackage{graphicx}
\definecolor{labelkey}{rgb}{0,0.08,0.45}
\definecolor{refkey}{rgb}{0,0.6,0.0}
\definecolor{Brown}{rgb}{0.45,0.0,0.05}
\definecolor{lime}{rgb}{0.00,0.8,0.0}
\definecolor{lblue}{rgb}{0.5,0.5,0.99}
\usepackage{pstricks-add}
\usepackage{amssymb}
\usepackage{fancyhdr}
\usepackage{tikz}
\usepackage{enumerate}
\usepackage[margin=1.1in]{geometry}
\usetikzlibrary{arrows}
\usepackage{hyperref}
\usepackage{esint}
\definecolor{labelkey}{rgb}{0.6,0.6,0.6}
\definecolor{refkey}{rgb}{0,0.6,0.0}
\def\disp{\displaystyle}
\def\e{\epsilon}

\def\Tilde{\widetilde}
\def\tilde{\widetilde}
\def\tx{\widetilde{x}}
\def\tu{\widetilde{u}}

\def\({\left(}
\def\){\right)}
\def\[{\left[}
\def\]{\right]}
\def\n{\left \|}
\def\en{\right \|}

\def\ox{\bar{x}}
\def\oy{\bar{y}}

\def\ou{\bar{u}}

\def\gg{\gamma}

\def\tilde{\widetilde}
\def\Tilde{\widetilde}

\def\la{\langle}
\def\ra{\rangle}

\def\h{\hfill\Box}
\def\R{\mathbb{R}}
\def\N{\mathbb{N}}

\def\co{\mbox{\rm co}\,}

\def\proj{\mbox{\rm proj}\,}

\def\gg{\gamma}

\def\ph{\varphi}

\def\N{I\!\!N}

\begin{document}

\def\proj{\text{proj}}
\author{Boris Mordukhovich\thanks{Department of Mathematics, Wayne State University, Detroit, MI 48201, USA; e-mail: aa1086@wayne.edu. Research of this author was partly supported by the US National Science Foundation under grant DMS-2204519, by the Australian Research Council under Discovery Project DP-190100555, and by Project~111 of China under grant D21024.}\;\;\;\;Dao Nguyen\thanks{Department of Mathematics and Statistics, San Diego State University, San Diego, CA, USA; email: dnguyen28@sdsu.edu}\;\;\;\;Trang Nguyen\thanks{Department of Mathematics, Wayne State University, Detroit, MI 48201, USA; e-mail: daitrang.nguyen@wayne.edu. Research of this author was partly supported by the US National Science Foundation under grant DMS-2204519.}\\Norma Ortiz-Robinson\thanks{Department of Mathematics, Grand Valley State  University, Allendale, MI 49401, USA; e-mail: ortizron@gmail.gvsu.edu.}\;\;\;\;\;Vinicio R\'ios\thanks{Department of Mathematics, Louisiana State University, Baton Rouge, LA 70803, USA; e-mail: vrios3@lsu.edu.}}

\title{Well-Posedness and Stability of Discrete Approximations for Controlled Sweeping Processes with Time Delay}
\maketitle
\thispagestyle{empty}
\vspace*{-0.2in}
\begin{center}
{\bf DEDICATED TO THE MEMORY OF ANDREA BACCIOTTI}
\end{center}

{\bf Abstract.} This paper addresses, for the first time in the literature, optimal control problems for dynamic systems governed by a novel class of sweeping processes with time delay. We establish well-posedness of such processes, in the sense of the existence and uniqueness of feasible trajectories corresponding to feasible controls under fairly unrestrictive assumptions. Then we construct a well-posed family of discrete approximations and find efficient conditions under the discretized time-delayed sweeping process exhibits stability with respect to strong convergence of feasible and optimal solutions. This creates a bridge between optimization of continuous-time and discrete-time sweeping control systems and justifies the effective use of discrete approximations in deriving optimality conditions and numerical techniques to solve the original time-delayed sweeping control problems via discrete approximations.



\section{Introduction}\label{intro}
Well-posedness and stability are among the most important issues in control theory, optimization and their applications. These and related aspects of optimal control in various settings have been largely addressed in the work by Andrea Bacciotti; see, e.g., his recent book \cite{A} and numerous publications including, in particular, differential inclusions, discontinuous dynamics, etc.

In this paper, we consider a new class of control systems the dynamics of which is governed by discontinuous differential inclusions introduced by Jean Jacques Moreau in the 1970s under the name of {\em sweeping processes}; see \cite{moreau} and the references therein. Moreau and his followers studied {\em uncontrolled} sweeping processes governed by ordinary differential inclusions of the form
\begin{equation}\label{sw}
-\dot x(t)\in N_{C(t)}\big(x(t)\big)\;\mbox{ a.e. }\;t\in[0,T],\quad x(0)=x_0\in C(0),
\end{equation}
via the normal cone to a continuously moving convex (or slightly nonconvex) set. Over the years, many researchers have been involved in developing sweeping process theory and its various applications to practical models in mechanics, hysteresis, traffic equilibria, etc. The reader is referred to the excellent, rather recent survey \cite{BT} with the vast bibliography therein. 

Optimal control problems for sweeping processes have been formulated and investigated much later; see \cite{ET} for the study of existence and relaxation issues, and \cite{BK,chhm12} for deriving necessary optimality conditions in different control models. It has been realized that optimal control problems for sweeping processes are highly challenging due to some phenomena that have never been exposed in control theory. They include the intrinsic discontinuity and unboundedness of the vector field in \eqref{sw}, the unavoidable presence of pointwise state and irregular mixed constraints on controls and trajectories, etc. 

Recent years have witnessed a strong interest in the study and applications of various optimal control problems for sweeping processes of type \eqref{sw} with control functions entering moving sets and additive perturbations of the dynamics. Among major contributions, we mention \cite{cm19,chhm16,cmn18,cmnn24,pfs,hp,vera}.

To the best of our knowledge, sweeping control problems with time delay have never been considered in the literature. We initiate such a study in this paper and will continue it in our further research. The {\em sweeping process with time-delayed controlled perturbations} considered below is formulated as follows:
\begin{equation*}
\begin{array}{ll}
-\dot{x}(t)&\in N_{C(t)}\big(x(t)\big)+g\big(t,x(t),x(t-\delta(t)),u(t)\big)\;\mbox{ a.e. }\;t\in I,\\
\;\;x(t)&\in C(t)\;\mbox{ for }\;t\in I,\quad x(s)=\phi(s)\;\mbox{ for }\;s\in[-\Delta,t_0],\\ \tag{$\mathcal{P}$}
\;\;u(t)&\in U\;\mbox{ a.e. }\;t\in I,
\end{array}
\end{equation*}
where the perturbation $g$ is a single-valued mapping, $\delta(\cdot)$ is a prescribed delay function that is nonnegative on the interval $I:=[t_0,T]$ with $T>t_0\ge 0$, $\phi(\cdot)$ is a given history defined on $[-\Delta,t_0]$ with $\Delta>0$, and the control set $U\subset\mathbb{R}^m$ is nonempty and compact. For every $u(\cdot)\in L^1(I,\mathbb{R}^m)$, a solution/trajectory to ($\mathcal{P}$) is a mapping $x:[-\Delta,T]\rightarrow\mathbb{R}^n$ that is absolutely continuous on $[t_0,T]$ satisfying the pointwise state constraint on $[0,T]$ and the initial condition on $[-\Delta,t_0]$ given in the second line above. 

Postponing the formulation of the precise assumption on the perturbation mapping $g$, we now clarify the meaning of the normal cone to $C(t)$ and the standing assumptions on the moving set/set-valued mapping $C(\cdot)$. The finite-dimensional notions and results of variational analysis and generalized differentiation used in what follows can be found in the standard texts \cite{CLSW,m18,rw}. 

The {\em proximal normal cone} to a closed set $S$ at $x\in S\subset\mathbb{R}^n$ is defined by
\begin{align*}
\tag{$\mathcal{NC}$}
N^P_S(x):=\big\{v\in\mathbb{R}^n\;\big|\;\langle v,y-x\rangle\le \sigma_v\|y-x\|^2\;\textup{ for all }\;y\in S\;\textup{ and some }\;\sigma_v\ge 0\big\}.
\end{align*}
Given $r>0$, a set $S\in\mathbb{R}^n$ is called (uniformly) $r$-{\em prox-regular} if for all $x,y\in S$ we have
$$
\langle \zeta, y-x\rangle\le \frac{1}{2r}\|y-x\|^2\;\mbox{ whenever }\;\zeta\in N_{S}^P(x).
$$
It has been well-recognized in the literature that for $r$-prox-regular sets, the major normal cones of variational analysis (Clarke, Fr\'echet, Mordukhovich) coincide at the point in question and reduce to the proximal normal cone ($\mathcal{NC}$). Therefore, the same normal cone notation $N_S$ can be used for such sets.

Taking this into account, we now list the standing assumptions of the moving set $C(t)$.\par
\noindent
{\bf(C1)} For each $t\in I$, the set $C(t)$ is nonempty and closed.\par
\noindent
{\bf(C2)} The mapping $C(\cdot)$ is absolutely continuous, which means that there exists a real-valued absolutely continuous function $v(\cdot)$ such that for all $y\in\mathbb{R}^n$ and all $s,t\in I$, we have \par
\noindent
$$
|d_{C(t)}(y)-d_{C(s)}(y)|\le |v(t)-v(s)|,
$$
where $d_S(\cdot):=\inf_{y\in S}\,\|y-\cdot\|$ is the distance function associated with the set $S$. \par
\noindent
{\bf(C3)} For each $t\in I$, the set $C(t)$ is $r$-prox regular.

Having (C3) in mind, recall the well-known relationship
\begin{equation}\label{nonsmooth}
N_S(x)\cap\mathbb{B}=\partial d_S(x)
\end{equation}
between the normal cone to a uniformly prox-regular set and the corresponding subdifferential of the distance function, which is used below for the case where $S=C(t)$.

Let us mention that control problems for differential inclusions with time delays of various types were studied in \cite[Chapter~7]{m06} and the references therein, but for {\em Lipschitzian} delay-differential inclusions and the like, where the Lipschitz continuity was absolutely crucial. Since the delay-differential inclusion in ($\mathcal{P}$) is highly non-Lipschitzian, the study of ($\mathcal{P}$) requires developing a significantly different machinery of variational analysis, which is done below.

 The structure of the paper is outlined as follows. In Section~2, we establish the {\em existence and uniqueness} of {\em feasible trajectories} $x(\cdot)$ for the perturbed time-delayed sweeping control problem ($\mathcal{P}$) corresponding to any feasible control function $u(\cdot)$. Moreover, an {\em explicit bound estimate} of $\|x(\cdot)\|$ in the $L^\infty$-norm is obtained in terms of the given problem data and the imposed assumptions.

 In Section~3, we consider the time-delayed sweeping control system in ($\mathcal{P}$) with a polyhedral set $C$ and construct a well-posed sequence of {\em discrete/finite-difference approximations} of this system in such a way that any feasible pair $(x(\cdot),u(\cdot))$ can be {\em strongly approximated} in the $W^{1,2}\times L^2$-norm by feasible solutions to the discrete-time controlled sweeping processes.

 In Section~4, we formulate a Mayer-type {\em optimal control problem} governed by time-delayed perturbed sweeping processes, introduce an appropriate notion of {\em local minimizers} for this problem, and construct a sequence of discrete-time optimization systems for which optimal solutions exist and {\em $W^{1,2}\times L^2$-strongly approximate} the designated local minimizer of the original sweeping control problem. This result serves as a critical link between continuous-time sweeping control problems with delay and their discrete approximations. In the concluding Section~5, we discuss some directions of our future research.
 
\section{Well-Posedness of Time-Delayed Sweeping Processes}
\setcounter{equation}{0}

In this section, we consider the time-delayed constrained sweeping process in ($\mathcal{P}$) for any fixed control $u(\cdot)$ and rewrite ($\mathcal{P}$) in the form
\begin{align*}
-\dot{x}(t)&\in N_{C(t)}\big(x(t)\big)+f\big(t,x(t),x(t-\delta(t))\big)\;\mbox{ a.e. }\;t\in I,\\
x(t)&\in C(t),\;t\in I,\tag{$\mathcal{AP}$}\\
x(s)&=\phi(s),\;s\in[-\Delta,t_0],   
\end{align*}
where $f(\cdot)=g(\cdot,u)$. In addition to the assumptions (C1)--(C3) on the moving set, let us list the {\em standing assumptions} on the other data in ({$\mathcal{AP}$):

{\bf(SH1)} The mapping $f$ is measurable in $t$, and for every $\eta>0$ there is $k_{\eta}(\cdot)\in L^1(I,\mathbb{R})$ such that
$$
\|f(t,x_1,y_1)-f(t,x_2,y_2)\|\le k_{\eta}(t)\{\|x_1 - x_2\| + \|y_1 - y_2\|\}\;\textup{ for all }\;t\in I
$$
whenever $(x_1,y_1), (x_2,y_2)\in \mathbb{B}[0,\eta]\times\mathbb{B}[0,\eta]$, where $\mathbb{B}(0,\eta)$ denotes the closed ball in $\mathbb{R}^n$ centred at the origin with radius $\eta>0$. 

{\bf(SH2)} There is $\beta(\cdot)\in L^1(I,\mathbb{R})$ such that we have the linear growth condition
\begin{equation*}
\left\|f(t,x,y)\right\|\le \beta(t)(1+\left\|x\right\|+\left\|y\right\|)\;\textup{ for all }\;(t,x,y)\in I\times C\times \mathbb{R}^n
\end{equation*}
with the notation $C:=\bigcup_{s\in I} C(s)$.

{{\bf(SH3)} The {\em delay function} $\delta:I\rightarrow [0,\infty)$ is Lipschitz continuous, nonincreasing, and satisfies the condition $0<\Delta:=\max\{\delta(t)\;|\;t\in I\}<\infty$.}

{\bf(SH4)} The {\em history function} $\phi:[-\Delta,t_0]\rightarrow\mathbb{R}^n$ is Lipschitz continuous.

The result below is taken from \cite[Proposition~1]{ET}.

\begin{lema} \label{simpleper}Let the assumptions in {\rm(C1)--(C3)} and {\rm(SH4)} be fulfilled. Then for every $h(\cdot)\in L^1(I,\mathbb{R}^n)$, the undelayed sweeping process 
\begin{align*}
-\dot{x}(t)&\in N_{C(t)}\big(x(t)\big)+h(t)\;\mbox{ a.e. }\;t\in I,\\
x(t)&\in C(t),\;t\in I,\\
x(s)&=\phi(s),\;s\in[-\Delta,t_0],   
\end{align*}
admits a unique absolutely continuous solution $x(\cdot)$ on $I$ satisfying the estimate
$$
\|\dot x(t) + h(t)\| \le \|h(t)\|+|\dot v(t)|\;\textup{ a.e. }\;t\in I.
$$
\end{lema}

Now we are ready to establish the following {\em well-posedness} theorem for the time-delayed sweeping process in ($\mathcal{AP}$) ensuring the existence and uniqueness of solutions together with an explicit bound estimate.

\begin{theor}\label{existence} Under the fulfillment of all the assumptions in {\rm(C1)--(C3)} and {\rm(SH1)--(SH3)}, there exists a unique absolutely continuous solution $x(\cdot)$ to ($\mathcal{AP}$) with the $L^{\infty}$-norm estimate $\|x(\cdot)\|_\infty\le l$, where
$$
l:=\|x_0\|+\exp\bigg\{4\int_{T_0}^T\beta(s)\,ds\bigg\}\int_{T_0}^T2\beta(s) (1+\|\phi\|_\infty+2\|x_0\|+|\dot{v}(s)|)\,ds.
$$
\end{theor}
{\bf Proof}. Suppose for convenience that
\begin{equation}\label{beta}
\int_{t_0}^T \beta(s)\,ds \le \frac{1}{8},
\end{equation}
which does not restrict the generality, and split the proof into the four major steps.

{\bf Step~1:} {\em Construction of approximate solutions}. At this step, we construct a  
sequence of functions $\{x_k(\cdot)\}$ in $\textup{C}(I,\mathbb{R}^n)$, which can be viewed as a sequence of approximating solutions to ($\mathcal{AP}$). To proceed, fix any $k\in\mathbb{N}$ and define the {\em uniform partition} $t_i^k:=t_0+i\frac{T-t_0}{k}$, $i=0,1,\ldots,k$, of the interval $[t_0,T]$. Then consider the sweeping process
\begin{align*}
-\dot{x}(t)&\in N_{C(t)}\big(x(t)\big)+f\big(t,\phi(t_0),\phi(t_0-\delta(t_0))\big)\;\mbox{ a.e. }\;t\in [t_0,t_1],\\
x(t)&\in C(t),\;t\in[t_0,t_1],\\
x(s)&=\phi(s),\,s\in[-\Delta,t_0].   
\end{align*}
Letting $h_0(t):=f(t,\phi(t_0),\phi(t_0-\delta(t_0)))$ allows us applying 
Lemma~\ref{simpleper} to obtain a unique absolutely continuous solution $x_0^k(\cdot):[t_0,t_1^k]\rightarrow\mathbb{R}^n$ of the above system that satisfies the estimate
$$
\|\dot x_0^k(t) + h_0(t)\|\le\|h_0(t)\|+|\dot v(t)|\;\textup{ a.e. }\;
t\in [t_0,t_1^k].
$$
Recursively, for each $i=1,\ldots,k-1$ let $p$ be such a number that $t_i^k-\delta(t_i^k)\in [t_{i-p}^k,t_{i-p+1}^k]$, and let $h_{i}(t):=f(t,x_{i-1}^k(t_i^k),x_{i-p}^k(t_i^k-\delta(t_i^k)))$. Applying again Lemma~\ref{simpleper} brings us to the system
\begin{align*}
-\dot{x}(t)&\in N_{C(t)}\big(x(t)\big)+h_i(t)\;\mbox{ a.e. }\;t\in[t_i^k,t_{i+1}^k],\\
x(t)&\in C(t),\;t\in [t_i^k,t_{i+1}^k],\\
x(s)&=\phi_i^k(s),\;s\in[-\Delta,t_i^k],   
\end{align*}
with $\phi_i^k(s)=\phi(s)$ for $s\in[-\Delta,t_0]$ and $\phi_i^k(s)=x_{r}^k(s)$ for $s\in[t_{r}^k,t_{r+1}^k]$ as $r=0,\ldots,i-1$. In this way, we get a unique absolutely continuous solution $x_i^k(\cdot):[t_i^k,t_{i+1}^k]\rightarrow\mathbb{R}^n$ to the latter system that satisfies
$$
\|\dot x_i^k(t) + h_i(t)\| \le \|h_i(t)\|+|\dot v(t)|\;\textup{ a.e. }\;t\in [t_i^k,t_{i+1}^k].
$$
Observe that $x_i^k(t_i^k)=x_{i-1}^k(t_i^k)$ for $i=1,\ldots,k$. For the consecutive arcs $x_{i-1}^k(\cdot)$ and $x_{i}^k(\cdot)$,  define
\begin{align*}
x_k(t):&=x_i^k(t)\;\textup{ if }\;t\in [t_i^k,t_{i+1}^k]\;\textup{ for some }\;i=0,1,\ldots,k-1,\\
x_k(t):&=\phi(t)\;\textup{ if }\;t\in [-\Delta,t_0],\;\mbox{ whenever }\;k\in\mathbb{N}.
\end{align*}
It is clear that $\{x_k(\cdot)\}$ is a sequence of arcs that are continuous on $[-\Delta,t_0]$ and absolutely continuous on $[t_0,T]$. Considering $\theta_k(t_0):=t_0$ and $\theta_k(t):=t_i^k$ if $t\in (t_i^k,t_{i+1}^k]$ for $i=0,\ldots,k-1$ leads us to
\begin{equation}\tag{a}
x_{i-1}^k(t_i^k)=x_i^k(t_i^k)=x_k(t_i^k)=x_k(\theta_k(t))\;\textup{ for all }\;t\in (t_i^k,t_{i+1}^k],
\end{equation}
\begin{equation}\tag{b}
x_{i-p}^k(t_i^k-\delta(t_i^k))=x_k(t_i^k-\delta(t_i^k))=x_k(\theta_k(t)-\delta(\theta_k(t)))\;\textup{ for all }\;t\in (t_i^k,t_{i+1}^k].
\end{equation}
Therefore, each arc $x_k(\cdot)$ satisfies the system
\begin{equation}\label{dynamics}
\begin{array}{ll}
-\dot{x}_k(t)&\in N_{C(t)}\big(x_k(t)\big)+f\big(t,x_k(\theta_k(t)),x_k(\theta_k(t)-\delta(\theta_k(t)))\big)\;\textup{ a.e. }\;t\in I,\\
\;\; x_k(t)&\in C(t),\;t\in I,\\
\;\;x_k(s)&=\phi(s),\;s\in[-\Delta,t_0],
\end{array}
\end{equation}
together with the following estimate valued for a.e. $t\in I$: 
\begin{equation}\label{bound1}
\big\|\dot x_k(t) + f\big(t,x_k(\theta_k(t)),x_k(\theta_k(t)-\delta(\theta_k(t)))\big)\big\| \le \big\|f\big(t,x_k(\theta_k(t)),x_k(\theta_k(t)-\delta(\theta_k(t)))\big)\big\|+|\dot v(t)|.
\end{equation}

{\bf Step~2:} {\em Convergence of approximating sequence}. At this step, we prove that the sequence $\{x_k(\cdot)\}$ converges uniformly to an absolutely continuous arc on $I$. To furnish this, let us show first that $\{\dot x_k(\cdot)\}$ is dominated by an integrable function in $I$. Employing \eqref{bound1} and the triangle inequality gives us for each $i\in\{0,1,\ldots,k-1\}$ that
\begin{equation}\label{derbound}
\begin{array}{ll}
 \|\dot x_k(t)\|&\le\big\|\dot x_k(t)+ f\big(t,x_k(t_i^k),x_k(t_i^k-\delta(t_i^k))\big)\big\|+\big\|f\big(t,x_k(t_i^k),x_k(t_i^k-\delta(t_i^k))\big)\big\|\\[1ex]
 &\le 2\big\|f\big(t,x_k(t_i^k),x_k(t_i^k-\delta(t_i^k))\big)\big\|+|\dot v(t)|\;\;\textup{a.e. }\;t\in [t_i^k,t_{i+1}^k].
 \end{array}
\end{equation}
Since $x_k(\cdot)$ is absolutely continuous on $I$, we have
$$
x_k(t_{i+1}^k)=x_k(t_i^k)+\int_{t_i^k}^{t_{i+1}^k}\,\dot x_k(t)\,dt
$$
and thus arrive at the estimates
\begin{align*}
\|x_k(t_{i+1}^k)\|&\le \|x_k(t_i^k)\|+ \int_{t_i^k}^{t_{i+1}^k}\|\dot x_k(t)\|\,dt\\
&\le \|x_k(t_i^k)\|+ 2\int_{t_i^k}^{t_{i+1}^k}\big\|f\big(t,x_k(t_i^k),x_k(t_i^k-\delta(t_i^k))\big)\big\|\,dt+\int_{t_i^k}^{t_{i+1}^k}|\dot v(t)|\,dt.
\end{align*}
Iterating on $\|x_k(t_i^k)\|$ readily brings us to
$$
\|x_k(t_{i+1}^k)\|\le\|x_0\| +2\sum_{r=0}^i\int_{t_r^k}^{t_{r+1}^k} \big\|f\big(t,x_k(t_r^k),x_k(t_r^k-\delta(t_r^k))\big)\big\|\,dt+\sum_{r=0}^i\int_{t_r^k}^{t_{r+1}^k}|\dot v(t)|\,dt.
$$
Recall that for every $r\in\{0,1,\ldots,k\}$ we have $(x_k(t_r^k),x_k(t_r^k-\delta(t_r^k)))\in C\times\mathbb{R}^n$. Then the linear growth condition in (SH2) ensures that
$$
\big\|f\big(t,x_k(t_r^k),x_k(t_r^k-\delta(t_r^k))\big)\big\|\le \beta(t)(1+\|x_k(t_r^k)\|+\|x_k(t_r^k-\delta(t_r^k))\|)\;\textup{ for all }\;t\in I.
$$
Getting all this together, we obtain
\begin{equation}\label{boundk}
\|x_k(t_{i+1}^k)\|\le\|x_0\| +2\sum_{r=0}^i(1+\|x_k(t_r^k)\|+\|x_k(t_r^k-\delta(t_r^k))\|)\int_{t_r^k}^{t_{r+1}^k}\beta(t)\,dt+\sum_{r=0}^i\int_{t_r^k}^{t_{r+1}^k}|\dot v(t)|\,dt.
\end{equation}
If $t_r^k-\delta(t_r^k)<t_0$, then $\|x_k(t_r^k-\delta(t_r^k))\|=\|\phi(t_r^k-\delta(t_r^k))\|\le \|\phi(\cdot)\|_{\infty}$. On the other hand, it follows from $t_r^k-\delta(t_r^k)\ge t_0$ and the absolute continuity of $x_k(\cdot)$ that
\begin{align*}
\|x_k(t_r^k-\delta(t_r^k))\|&\le \|x_k(t_r^k)\|+\int_{t_r^k}^{t_r^k-\delta(t_r^k)}\|\dot x_k(t)\|dt\\
&=\|x_k(t_r^k)\|-\int_{t_r^k-\delta(t_r^k)}^{t_r^k}\|\dot x_k(t)\|dt\\
&\le \|x_k(t_r^k)\|.
\end{align*}
Regardless the value of $t_r^k-\delta(t_r^k)$, we get the estimate
\begin{equation}\label{bound}
\|x_k(t_r^k-\delta(t_r^k))\|\le\|\phi(\cdot)\|_{\infty}+\|x_k(t_r^k)\|.
\end{equation}
According to \eqref{boundk}, the latter implies by \eqref{beta} that
\begin{align*}
\|x_k(t_{i+1}^k)\|&\le\|x_0\| +2\sum_{r=0}^i(1+\|\phi(\cdot)\|_{\infty}+2\|x_k(t_r^k)\|)\int_{t_r^k}^{t_{r+1}^k}\beta(t)\,dt+\sum_{r=0}^i\int_{t_r^k}^{t_{r+1}^k}|\dot v(t)|\,dt\\
&\le\|x_0\| +2(1+\|\phi(\cdot)\|_{\infty}+2\hspace{-0.05cm}\max_{0\le r\le k}\|x_k(t_r^k)\|)\sum_{r=0}^i\int_{t_r^k}^{t_{r+1}^k}\beta(t)\,dt+\sum_{r=0}^i\int_{t_r^k}^{t_{r+1}^k}|\dot v(t)|\,dt\\
&\le \|x_0\| + 2(1 + \|\phi(\cdot)\|_{\infty}+2\hspace{-0.05cm}\max_{0\le r\le k}\|x_k(t_r^k)\|)\int_{t_0}^{T}\beta(t)\,dt+\int_{t_0}^{T}|\dot v(t)|\,dt\\
&\le \|x_0\| + \frac{1}{4}(1+\|\phi(\cdot)\|_{\infty}) + \frac{1}{2}\max_{0\le r\le k}\|x_k(t_r^k)\|+\int_{t_0}^{T}|\dot v(t)|\,dt,
\end{align*}
which holds for any $i\in\{0,1,\ldots,k-1\}$. Therefore,
$$
\max_{0\le r\le k}\|x_k(t_r^k)\|\le \|x_0\| + \frac{1}{4}(1 + \|\phi(\cdot)\|_{\infty})+ \frac{1}{2}\max_{0\le r\le k}\|x_k(t_r^k)\|+\int_{t_0}^{T}|\dot v(t)|\,dt,\;\mbox{ i.e.,}
$$
\begin{equation}\label{xmax}
\max_{0\le r\le k}\|x_k(t_r^k)\|\le 2\bigg(\|x_0\| + \frac{1}{4}(1+\|\phi(\cdot)\|_{\infty})+\int_{t_0}^{T}|\dot v(t)|\,dt\bigg)=:M.
\end{equation}
Using (\ref{derbound}), the linear growth condition in (SH2), (\ref{bound}), and \eqref{xmax} tells us that
\begin{equation}\label{boundder}
\begin{array}{ll}
\|\dot x_k(t)\|&\le 2\|f(t,x_k(t_i^k),x_k(t_i^k-\delta(t_i^k)))\|+|\dot v(t)|\\[1ex]
&\le2(1+\|\phi(\cdot)\|_{\infty}+2M)\beta(t)+|\dot v(t)|\;\textup{ a.e. }\;t\in I. 
\end{array}
\end{equation}
In what follows, we use estimate (\ref{boundder}) to show that $\{x_k(\cdot)\}$ is a Cauchy sequence with respect to the uniform topology. Employing the conditions in (a), (b), (\ref{bound1}), and (\ref{boundder}) leads us to the estimate
\begin{equation}\label{twovelocities}
\big\|\dot x_k(t) + f\big(t,x_k(\theta_k(t)),x_k(\theta_k(t)-\delta(\theta_k(t)))\big)\big\|\le (1+\|\phi(\cdot)\|_{\infty}+2M)\beta(t)+|\dot v(t)|\;\textup{ a.e. }\;t\in I.   
\end{equation}
Define further $\alpha(t)=(1+\|\phi(\cdot)\|_{\infty}+2M)\beta(t)+|\dot v(t)|$ and deduce from (\ref{dynamics}) for any $k,l\in\mathbb{N}$ that
$$
-\frac{1}{\alpha(t)}\big[\dot x_k(t) + f\big(t,x_k(\theta_k(t)),x_k(\theta_k(t)-\delta(\theta_k(t)))\big)\big]\in N_{C(t)}\big(x_k(t)\big)\;\textup{ a.e. }\;t\in I,
$$
$$
-\frac{1}{\alpha(t)}\big[\dot x_l(t) + f\big(t,x_l(\theta_l(t)),x_l(\theta_l(t)-\delta(\theta_l(t)))\big)\big]\in N_{C(t)}\big(x_l(t)\big)\;\textup{ a.e. }\;t\in I.
$$
Since $C(t)$ is $r$-prox-regular, it follows from the latter inclusions that
\begin{align*}
\big\langle\dot x_k(t) &+ f\big(t,x_k(\theta_k(t)),x_k(\theta_k(t)-\delta(\theta_k(t)))\big)-\dot x_l(t)-f\big(t,x_l(\theta_l(t)),x_l(\theta_l(t)-\delta(\theta_l(t)))\big),x_k(t)-x_l(t)\big\rangle\\
&=\big\langle-\big[\dot x_k(t) + f\big(t,x_k(\theta_k(t)),x_k(\theta_k(t)-\delta(\theta_k(t))\big)\big],x_l(t)-x_k(t)\big\rangle\\
&+\big\langle-\big[\dot x_l(t) + f\big(t,x_l(\theta_l(t)),x_l(\theta_l(t)-\delta(\theta_l(t)))\big)\big],x_k(t)-x_l(t)\big\rangle\\
&\le \frac{\alpha(t)}{r}\big\|x_l(t)-x_k(t)\big\|^2,
\end{align*}
from which we clearly get that 
\begin{align*}
\big\langle\dot x_k(t)&-\dot x_l(t),x_k(t)-x_l(t)\big\rangle\\
&\le\big\langle  f\big(t,x_k(\theta_k(t)),x_k(\theta_k(t)-\delta(\theta_k(t)))\big)-f\big(t,x_l(\theta_l(t)),x_l(\theta_l(t)-\delta(\theta_l(t)))\big),x_l(t)-x_k(t)
\big\rangle\\
&+\frac{\alpha(t)}{r}\big\|x_l(t)-x_k(t)\big\|^2. 
\end{align*}
Therefore, the above relationships bring us to the inequality
\begin{equation}\label{quadratic}
\begin{array}{ll}
&\disp\frac{1}{2}\frac{d}{dt}\|x_k(t)-x_l(t)\|^2=\langle\dot x_k(t)-\dot x_l(t),x_k(t)-x_l(t)\big\rangle\\[1.5ex]
&\le\big\langle f\big(t,x_k(\theta_k(t)),x_k(\theta_k(t)-\delta(\theta_k(t)))\big)-f\big(t,x_k(t),x_k(t-\delta(t))\big),x_l(t)-x_k(t)\big\rangle\\[1ex]
&+\big\langle f\big(t,x_k(t),x_k(t-\delta(t)))-f(t,x_l(t),x_l(t-\delta(t))\big),x_l(t)-x_k(t)\big\rangle\\[1ex]
&+\big\langle f\big(t,x_l(t),x_l(t-\delta(t))\big)-f\big(t,x_l(\theta_l(t)),x_l(\theta_l(t)-\delta(\theta_l(t)))\big),x_l(t)-x_k(t)\big\rangle\\[1ex]
&+\disp\frac{\alpha(t)}{r}\|x_l(t)-x_k(t)\|^2.
\end{array}
\end{equation}
It follows from (\ref{boundder}) that we have
\begin{align*}
\disp\|x_k(t)\|\le \|x_0\|+\int_{t_0}^t \|\dot x_k(s)\|\,ds&\le \|x_0\|+2(1+\|\phi(\cdot)\|_{\infty}+2M)\int_{t_0}^t \beta(s)\,ds+\int_{t_0}^t |\dot v(s)|\,ds\\
&\le \|x_0\|+\frac{1}{4}(1+\|\phi(\cdot)\|_{\infty}+2M)+\int_{t_0}^T |\dot v(s)|\,dt=:\eta>0,
\end{align*}
which means that $x_k(t)\in\mathbb{B}[0,\eta]$ for all $t\in I$ and all $k\in\mathbb{N}$. Hence $\|x_k(t)-x_l(t)\|\le 2\eta$ whenever $t\in I$. By (SH1), there exists $k(\cdot)\in L^1(I,\mathbb{R})$ such that $f(t,\cdot,\cdot)$ is $k(t)$-Lipschitz on $\mathbb{B}[0,\eta]\times\mathbb{B}[0,\eta]$. From (\ref{quadratic}) we therefore deduce the following estimate valid for a.e. $t\in I$: 
\begin{align*}
&\frac{1}{2}\frac{d}{dt}\|x_k(t)-x_l(t)\|^2\\
&\le k(t)\{\|x_k(\theta_k(t))-x_k(t)\|+\|x_k(\theta_k(t)-\delta(\theta_k(t)))-x_k(t-\delta(t))\|\}\|x_l(t)-x_k(t)\|\\
&+k(t)\{\|x_k(t)-x_l(t)\|+\|x_k(t-\delta(t))-x_l(t-\delta(t))\|\}\|x_l(t)-x_k(t)\|\\
&+k(t)\{\|x_l(t)-x_l(\theta_l(t))\|+\|x_l(t-\delta(t))-x_l(\theta_l(t)-\delta(\theta_l(t)))\|\}\|x_l(t)-x_k(t)\|\\
&+\frac{\alpha(t)}{r}\|x_l(t)-x_k(t)\|^2.
\end{align*}
Collecting the like terms and making equivalent transformations yield

\begin{equation}\label{Estim}
\begin{array}{ll}
&\disp\frac{1}{2}\frac{d}{dt}\|x_k(t)-x_l(t)\|^2\\[1.5ex]
&\disp\le\bigg(k(t)+\frac{\alpha(t)}{r}\bigg)\|x_k(t)-x_l(t)\|^2+k(t)\|x_k(t-\delta(t))-x_l(t-\delta(t))\|\cdot\|x_k(t)-x_l(t)\|\\[1ex]
&+2\eta k(t)\{\|x_k(\theta_k(t))-x_k(t)\|+\|x_l(t)-x_l(\theta_l(t))\|\}\\[1ex]
&+2\eta k(t)\{\|x_k(\theta_k(t)-\delta(\theta_k(t)))-x_k(t-\delta(t))\|+\|x_l(t-\delta(t))-x_l(\theta_l(t)-\delta(\theta_l(t)))\|\}.
\end{array}
\end{equation}
\noindent
If $t-\delta(t)>0$, then for every $k,l\in\mathbb{N}$  we get from the absolute continuity of $x_k(\cdot)-x_l(\cdot)$ that
\begin{align*}
\|x_k(t-\delta(t))-x_l(t-\delta(t))\|&\le\|x_k(t)-x_l(t)\|+\int_{t}^{t-\delta(t)}\|\dot x_k(s) - \dot x_l(s)\|\,ds\\
&=\|x_k(t)-x_l(t)\|-\int_{t-\delta(t)}^{t}\|\dot x_k(s) - \dot x_l(s)\|\,ds\\
&\le \|x_k(t)-x_l(t)\|.
\end{align*}
On the other hand, assuming $t-\delta(t)\le 0$ implies that $\|x_k(t-\delta(t))-x_l(t-\delta(t))\|=\|\phi(t-\delta(t))-\phi(t-\delta(t))\|=0.$ Therefore, we have for all $t\in I$ that
\begin{equation}\label{est1}
\|x_k(t-\delta(t))-x_l(t-\delta(t))\|\le\|x_k(t)-x_l(t)\|.    
\end{equation}
Defining $\gamma(t):=2(1+\|\phi(\cdot)\|_{\infty}+2M)\beta(t)+|\dot v(t)|$ on $I$, observe by (\ref{boundder}) that
\begin{equation}\label{estim2}
\|x_k(t)-x_k(\theta_k(t))\|\le\int_{\theta_k(t)}^t \|\dot x_k(s)\|\,ds\le\int_{\theta_k(t)}^t \gamma(s)\,ds\;\mbox{ whenever }\;t\in I,\;k\in\mathbb N.
\end{equation}
If $t-\delta(t)>t_0$, it follows from by continuity of $\delta(\cdot)$ and the 
convergence $\theta_k(t)\to t$ as $k\to\infty$ that $\theta_k(t)-\delta(\theta_k(t))\to t-\delta(t)$ and $\theta_k(t)-\delta(\theta_k(t))>t_0$ telling us that
\begin{equation*}
\|x_k(\theta_k(t)-\delta(\theta_k(t)))-x_k(t-\delta(t))\|\le \int_{\theta_k(t)-\delta(\theta_k(t))}^{t-\delta(t)} \|\dot x_k(s)\|\,ds\le\int_{\theta_k(t)-\delta(\theta_k(t))}^{t-\delta(t)}\gamma(s)\,ds.   
\end{equation*}
Similarly we get for $t-\delta(t)<t_0$ that $\theta_k(t)-\delta(\theta_k(t))\to t-\delta(t)$ and $\theta_k(t)-\delta(\theta_k(t))<t_0$ yielding
\begin{equation*}
\|x_k(\theta_k(t)-\delta(\theta_k(t)))-x_k(t-\delta(t))\|=\|\phi(\theta_k(t)-\delta(\theta_k(t)))-\phi(t-\delta(t))\|   
\end{equation*}
when $k$ is sufficiently large. Let us further define
\begin{equation}
\mathcal{D}_k(t):=
\begin{cases}\label{defD}
\disp\int_{\theta_k(t)-\delta(\theta_k(t))}^{t-\delta(t)}\gamma(s)\,ds  &\textup{ if }t-\delta(t)>t_0\\
\|\phi(\theta_k(t)-\delta(\theta_k(t)))-\phi(t-\delta(t))\| & \textup{ if }t-\delta(t)\le t_0
\end{cases}
\end{equation}
and observe that we have for large $k$ that
\begin{equation}\label{Destim}
|\mathcal{D}_k(t)|\le\max\bigg\{2\|\phi(\cdot)\|_{\infty},\, \int_{t_0}^T \gamma(s)\,ds 
\bigg\}=:\mathcal{M}, \;\;\lim_{k\to\infty} \mathcal{D}_k(t)=0\;\textup{ as }\;t\in I.  
\end{equation}
Plugging (\ref{est1})--(\ref{defD}) into (\ref{Estim}) gives us the estimate
\begin{align*}
\frac{1}{2}\frac{d}{dt}\|x_k(t)-x_l(t)\|^2&\le\bigg(2k(t)+\frac{\alpha(t)}{r}\bigg)\|x_k(t)-x_l(t)\|^2\\
&+2\eta k(t)\bigg\{\bigg(\int_{\theta_k(t)}^t \gamma(s)\,ds+\mathcal{D}_k(t)
\Big)+\Big(\int_{\theta_l(t)}^t \gamma(s)\,ds+\mathcal{D}_l(t)\Big)\bigg\}.
\end{align*}
Consider further the family of functions
$$
\mathcal{G}_{k,l}(t):=4\eta k(t)\bigg\{\bigg(\int_{\theta_k(t)}^t \gamma(s)\,ds+\mathcal{D}_k(t)\bigg)+\bigg(\int_{\theta_l(t)}^t \gamma(s)\,ds+\mathcal{D}_l(t)\bigg)\bigg\},\;\;t\in [t_0,T]
$$
for which we clearly get that
$$
\frac{d}{dt}\|x_k(t)-x_l(t)\|^2\le 4\bigg(k(t)+\frac{\alpha(t)}{2r}\bigg)\|x_k(t)-x_l(t)\|^2+\mathcal{G}_{k,l}(t)\;\textup{ a.e. }\;t\in I.
$$
Since $\gamma(\cdot)\in L^1(I,\mathbb{R})$ and $\theta_k(t)\to t$ for each $t\in I$, it follows from (\ref{Destim}) that
$$
\lim_{k,l\to \infty}\,\mathcal{G}_{k,l}(t)=0,\;\;|\mathcal{G}_{k,l}(t)|\le 8\eta\bigg(\mathcal{M}+\int_{t_0}^T\gamma(s)\,ds\bigg)k(t),\quad t\in I.
$$
We deduce from the Lebesgue dominated convergence theorem that
$$
\lim_{k,l\to \infty}\;\int_{t_0}^T \mathcal{G}_{k,l}(s)\,ds = 0.
$$
The above properties of $\mathcal{G}_{k,l}(\cdot)$ and the initial condition $\|x_k(t_0)-x_l(t_0)\|=0$ allow to apply \cite[Lemma~1]{ET} to guarantee the convergence
$$
\lim_{k,l\to \infty}\;\|x_k(\cdot)-x_l(\cdot)\|_{\infty}=0,
$$
which tells us that $\{x_k(\cdot)\}$ is a Cauchy sequence in $\textup{C}(I,\mathbb{R}^n)$ with respect to the uniform topology, and hence there exists $x(\cdot)\in \textup{C}(I,\mathbb{R}^n)$ such that $\lim_{k\to\infty}\|x_k(\cdot)-x(\cdot)\|_{\infty}=0$. Employing (\ref{boundder}) gives us a subsequence $\{x_{k_j}(\cdot)\}$ that converges weakly in $L^1(I,\mathbb{R}^n)$ to some $g(\cdot)\in L^1(I,\mathbb{R}^n)$. This allows us to write without relabeling that
$$
\int_{t_0}^t \dot x_{k}(s)\,ds \rightarrow \int_{t_0}^t g(s)\,ds\;\textup{ weakly in }\;\mathbb{R}^n\;\mbox{ whenever }\;t\in I.
$$
Therefore, for any $t\in I$ we have the equalities
\begin{align*}
x(t)&=\lim_{k\to\infty}\,x_k(t)\\
&=\lim_{k\to\infty}\,\bigg(x_0+\int_{t_0}^t\dot x_k(s)\,ds\bigg)\\
&=x_0+\lim_{k\to\infty}\,\int_{t_0}^t\dot x_k(s)\,ds\\
&=x_0+\int_{t_0}^t g(s)\,ds,
\end{align*}
which readily imply that $x(\cdot)$ is absolutely continuous on $I$, $\dot x(t)=g(t)$ a.e. $t\in I$, and $\dot x_k(\cdot)\to \dot x(\cdot)$ weakly in $L^1(I,\mathbb{R}^n)$ as $k\to\infty$.

{\bf Step~3:} {\em The limiting arc $x(\cdot)$ from Step~{\rm 2} is a unique solution to the sweeping process ($\mathcal{AP}$)}. First we verify that $x(\cdot)$ is a solution to ($\mathcal{AP}$). Since $\theta_k(t)\to t$ for any $t\in I$, the sequence $x_k(\cdot)$ converges uniformly to $x(\cdot)$, and  the functions $\delta(\cdot)$ and $\phi(\cdot)$ are continuous, we have that $x(\theta_k(t))\to x(t)$ and $x(\theta_k(t)-\delta(\theta_k(t)))\to x(t-\delta(t))$ for any $t\in I$. It follows from the continuity of $f(t,\cdot,\cdot)$ that
\begin{equation}\label{cont}
f\big(t,x(t),x(t-\delta(t))\big)= \lim_{k\to\infty}\,f\big(t,x_k(\theta_k(t)),x_k(\theta_k(t)-\delta(\theta_k(t)))\big).   
\end{equation}
Furthermore, we deduce from (a), (b), and (\ref{boundder}) that
\begin{equation}\label{boundfk}
 \|f(t,x_k(\theta_k(t)),x_k(\theta_k(t)-\delta(\theta_k(t))))\|\le  (1+\|\phi(\cdot)\|_{\infty}+2M)\beta(t),  
\end{equation}
which clearly ensures the estimate
\begin{equation}\label{boundf}
\begin{array}{ll}
\disp\|f(t,x(t),x(t-\delta(t)))\|&=\disp\lim_{k\to\infty}\,\|f(t,x_k(\theta_k(t)),x_k(\theta_k(t)-\delta(\theta_k(t))))\|\\
&\le\disp\lim_{k\to\infty}\,(1+\|\phi(\cdot)\|_{\infty}+2M)\beta(t)\\
&=(1+\|\phi(\cdot)\|_{\infty}+2M)\beta(t),\quad t\in I.
\end{array}
\end{equation}
Next we verify the fulfillment of the inclusion
$$
\dot x(t) + f\big(t,x(t),x(t-\delta(t))\big)\in N_C\big(x(t)\big)\;\mbox{ a.e. }\;t\in I.
$$
To furnish this, deduce from (\ref{cont})--(\ref{boundf}) and the weak convergence of $\{\dot x_k(\cdot)\}$ in $L^1(I,\mathbb{R}^n)$ that the sequence $\{\dot x_k(\cdot)+f(\cdot,x_k(\theta_k(\cdot)),x_k(\theta_k(\cdot)-\delta(\theta_k(\cdot))))\}$ is also weakly convergent in $L^1(I,\mathbb{R}^n)$. Then Mazur's weak closure theorem yields the existence of a sequence $\{\zeta_k(\cdot)\}$ that converges strongly in $L^1(I,\mathbb{R}^n)$ to $\dot x(\cdot)+f(\cdot,x(\cdot),x(\cdot-\delta(\cdot)))$ being such that 
$$
\zeta_k(t)\in\textup{co}\,\big\{\dot x_i(t)+f\big(t,x_i(\theta_i(t)),x_i(\theta_i(t)-\delta(\theta_i(t)))\big)\;\big|\;i\ge k\big\}
$$
for each $k$ and for all $t\in I$. Passing to a subsequence if necessary, we get without relabeling that
$$
\zeta_k(t)\rightarrow \dot x(t) + f(t,x(t),x(t-\delta(t)))\textup{ as }\;k\to\infty\;\textup{ a.e. }\;t\in I.
$$
Therefore, there exist finitely many functions $\lambda_i^k(t)\in[0,1]$, not more than $n+1$ for each $k$ by the Carath\'eodory theorem, with $\sum_{i\ge k}\lambda_i^k(t)=1$ and
\begin{equation}\label{convex}
 \zeta_k(t)=\sum_{i\ge k}\lambda_i^k(t)\big(\dot x_i(t)+f(t,x_i(\theta_i(t)),x_i(\theta_i(t)-\delta(\theta_i(t))))\big)
\end{equation}
for a.e. $t\in I$. Fix $j\in\mathbb{N}$ and observe that
$$
\zeta_k(t)\in\textup{co}\,\big\{\dot x_i(t)+f(t,x_i(\theta_i(t)),x_i(\theta_i(t)-\delta(\theta_i(t))))\;\big|\;i\ge j\big\}\;\textup{ whenever }\;k\ge j.
$$
Therefore, for a.e. $t\in I$ we have the conditions
\begin{align*}
\dot x(t) + f\big(t,x(t),x(t-\delta(t))\big)&=\lim_{k\to\infty}\,\zeta_k(t)\\
&\in\overline{\textup{co}}\,\big\{\dot x_i(t)+f(t,x_i(\theta_i(t)),x_i(\theta_i(t)-\delta(\theta_i(t))))\;\big|\;i\ge j\big\},\quad j\in\mathbb N,
\end{align*}
which imply in turn the inclusion
$$
\dot x(t) + f\big(t,x(t),x(t-\delta(t))\big)\in\bigcap_{k}\, \overline{\textup{co}}\,\big\{\dot x_i(t)+f\big(t,x_i(\theta_i(t)),x_i(\theta_i(t)-\delta(\theta_i(t)))\big)\;\big|\;i\ge k\big\}
$$
for a.e. $t\in I$. Picking $\xi\in\mathbb{R}^n$, we deduce from (\ref{convex}) that
\begin{equation}\label{hyper}
\begin{array}{ll}
\disp&\langle \xi,\dot x(t) + f(t,x(t),x(t-\delta(t)))\rangle\\
&=\disp\lim_{k\to\infty}\big\langle\xi,\sum_{i\ge k}\lambda_i^k(t)(\dot x_i(t)+f(t,x_i(\theta_i(t)),x_i(\theta_i(t)-\delta(\theta_i(t)))))\big\rangle\\ 
&=\disp\lim_{k\to\infty}\sum_{i\ge k}\lambda_i^k(t)\big\langle\xi,\dot x_i(t)+f(t,x_i(\theta_i(t)),x_i(\theta_i(t)-\delta(\theta_i(t))))\big\rangle\\
&\le\disp\lim_{k\to\infty}\sum_{i\ge k}\lambda_i^k(t)\sup_{r\ge k}\big\{\langle\xi,\dot x_r(t)+f(t,x_r(\theta_r(t)),x_r(\theta_r(t)-\delta(\theta_r(t))))\rangle\big\}\\
&=\disp\inf_{k}\,\sup_{r\ge k}\big\{\langle\xi,\dot x_r(t)+f(t,x_r(\theta_r(t)),x_r(\theta_r(t)-\delta(\theta_r(t))))\rangle\big\}\;\mbox{ a.e. }\;t\in I.
\end{array}
\end{equation}
It clearly follows from (\ref{twovelocities}) that 
\begin{equation}\label{al}
\frac{1}{\alpha(t)}\big(\dot x_k(t)+f(t,x_k(\theta_k(t)),x_k(\theta_k(t)-\delta(\theta_k(t)))\big)\in\mathbb{B},
\end{equation}
where $\mathbb B$ stands for the closed unit ball in $\mathbb{R}^n$. Since the sets $C(t)$ are uniformly prox-regular, combining
(\ref{nonsmooth}), (\ref{dynamics}), and \eqref{al} yields the relationship
$$
\frac{1}{\alpha(t)}\big(\dot x_k(t)+f(t,x_k(\theta_k(t)),x_k(\theta_k(t)-\delta(\theta_k(t)))\big)\in -N_{C(t)}(x_k(t))\cap\mathbb{B}=-\partial d_{C(t)}(x_k(t))
$$
for a.e. $t\in I$, and thus for such $t$ we get the inclusion
$$
\dot x_k(t)+f(t,x_k(\theta_k(t)),x_k(\theta_k(t)-\delta(\theta_k(t))))\in-\alpha(t)\hspace{0.01cm}\partial d_{C(t)}(x_k(t)).
$$
\begin{align*}
&\langle\xi,\dot x_r(t)+f(t,x_r(\theta_r(t)),x_r(\theta_r(t)-\delta(\theta_r(t))))\rangle\\
&\le\sup_{v\in -\alpha(t)\hspace{0.01cm}\partial d_{C(t)}(x_k(t))}\,\langle\xi,v\rangle\\
&=\sigma\big(-\alpha(t)\hspace{0.01cm}\partial d_{C(t)}(x_k(t)),\xi\big)\\
&=\alpha(t)\sigma\big(-\partial d_{C(t)}(x_k(t)),\xi\big),\quad\mbox{ a.e.},
\end{align*}
where $\sigma_S=\sigma(S)$ signifies the support function of the set in question. Using (\ref{hyper}) and the upper semicontinuity of the support function $\sigma(-\partial d_{C(t)}(\cdot),\xi)$ for the subdifferential of the distance function associated with uniformly prox-regular sets ensures the conditions 
\begin{align*}
\langle \xi,\dot x(t) + f(t,x(t),x(t-\delta(t)))\rangle & \le\alpha(t)\inf_{k}\,\sup_{r\ge k}\sigma(-\partial d_{C(t)}(x_k(t)),\xi)\\
&=\alpha(t)\limsup_{k}\sigma(-\partial d_{C(t)}(x_k(t)),\xi)\\
&\le\alpha(t)\sigma(-\partial d_{C(t)}(x(t)),\xi)
\end{align*}
for a.e. $t\in I$ and all $\xi\in\mathbb{R}^n$, which imply that for such t and $\xi$ we have
$$
\frac{1}{\alpha(t)}\big\langle\xi,\dot x(t) + f(t,x(t),x(t-\delta(t)))\big\rangle\le \sigma(-\partial d_{C(t)}(x(t)),\xi).
$$
Since $\partial d_{C(t)}(x(t))$ is closed and convex for uniformly regular sets, we get by (\ref{nonsmooth}) that
$$
\frac{1}{\alpha(t)}(\dot x(t) + f(t,x(t),x(t-\delta(t))))\in -\partial d_{C(t)}(x(t)),\xi)=-N_{C(t)}(x(t))\cap\mathbb{B},
$$
which tells us therefore that
$$
\dot x(t) + f(t,x(t),x(t-\delta(t))))\in -\alpha(t)\partial d_{C(t)}(x(t)),\xi)=-N_{C(t)}(x(t))\;\textup{ a.e. }\;t\in I.
$$
Taking into account that $x(t_0)=\lim_{k\to\infty}\,x_k(t_0)=\phi(0)$ justifies that $x(\cdot)$ is a solution to ($\mathcal{AP}$).

The uniqueness of solutions to ($\mathcal{AP}$) follows by the reduction arguments above from Lemma~\ref{simpleper}, being based on the {\em hypomonotonicity} of the normal cone mapping generated by $r$-prox-regular sets; cf.\ \cite{ET}. 

{\bf Step~4:} {\em Verifying the norm estimate}. Finally, we turn our attention to establishing the claimed estimate for this solution. It follows from Lemma~\ref{simpleper} and the linear growth condition that
\begin{align*}
\|\dot{x}(t)\|&\le 2\|f(t,x(t),x(t-\delta(t)))\|+|\dot{v}(t)|\\
&\le\beta(t)(1+\|x(t)\|+\|x(t-\delta(t))\|) +|\dot{v}(t)|.
\end{align*}
If $t-\delta(t)<0$, then $\|x(t-\delta(t))\|=\|\phi(t-\delta(t))\|\le\|\phi\|_\infty$, while for $t-\delta(t)\geq 0$ we have
\begin{align*}
\|x(t-\delta(t)\|&\le\|x(t)\|+\int_t^{t-\delta(t)}\|\dot{x}(s)\|\,ds\\
&=\|x(t)\|-\int_{t-\delta(t)}^t\|\dot{x}(s)\|\,ds\le\|x(t)\|.
\end{align*}
This readily implies that
\begin{equation*}
 \|\dot{x}(t)\|\leq 2 \beta(t)(1+\|\phi\|_\infty+2\|x(t)\|)+|\dot{v}(t)|,
\end{equation*}
and consequently, by using $\|x(t)\|\le\|x_0\|+\int_{T_0}^t \|\dot{x}(s)\|\,ds$, that
\begin{equation}
 \|\dot{x}(t)\|\leq 4\beta(t)\int_{T_0}^t \|\dot{x}(s)\|\,ds+ 2\beta(t)(1+\|\phi\|_\infty+2\|x_0\|)+|\dot{v}(t)|.
\end{equation}
Applying now Gronwall's lemma and remembering that $\beta(t)\geq 0$ give us
\begin{align*}
\|x(t)\|&\le\|x_0\|+\int_{T_0}^t 2\beta(s)(1+\|\phi\|_\infty+2\|x_0\|+|\dot{v}(s)|)\,\exp\bigg\{4\int_{s}^t\beta(\tau)\,d\tau\bigg\}\,ds\\
&\le\|x_0\|+\exp\bigg\{4\int_{T_0}^T\beta(s)\,ds\bigg\}\int_{T_0}^T2\beta(s) (1+\|\phi\|_\infty+2\|x_0\|+|\dot{v}(s)|)\,ds:=l,
\end{align*}
which justifies the claimed estimate and thus completes the proof of the theorem. $\h$

\section{Discrete Approximations of Sweeping  Processes with Delay}
\setcounter{equation}{0}

In this section, we formulate a Mayer-type optimal control problem for a controlled time-delayed sweeping process and construct a well-posed family of its discrete approximations. Our main goal here is to verify the possibility to approximate {\em any feasible} solution to the controlled sweeping process by a sequence of feasible solutions 
to its discrete-time counterparts. This is a principal result, which justifies the {\em stability} of controlled time-delayed sweeping processes with respect to discrete approximations. It definitely has a certain numerical flavor by showing that a complex  infinite-dimensional controlled dynamic systems can be eventually replaced by solving a finite-dimensional one, which admits a much broader arsenal of numerical methods. Moreover, the main result of this section plays a crucial role in establishing a strong approximation of a designated local minimizer of the time-delayed sweeping control problem established in the next section and in deriving necessary optimality conditions in our future research.

The time-delayed {\em sweeping optimal control problem} studied in this paper is formulated as follows:
\begin{equation}\tag{$\mathcal{DMP}$}
\begin{aligned}
\textup{minimize}\;\;J[x,u]&=\varphi(x(T))\\
\qquad\qquad\textup{subject to}\;-\dot x(t)&\in N_C\big(x(t)\big)+ g\big(t, x(t),x(t-\delta(t)),u(t)\big)\;\mbox{ a.e. }\;t\in[0,T]\\
x(s)&=\phi(s),\;s\in [-\Delta,0],\\
u(t)&\in U\;\mbox{ a.e. }\;t\in[0,T].
\end{aligned}
\end{equation}
where the convex polyhedral set $C$ is defined as
\begin{equation}\label{C}
C:=\bigcap_{j=1}^s C_j,\;\textup{ with }\;C_j:=\big\{x\in\mathbb{R}^n\;\big|\;\langle x_*^j,x\rangle\le c_j\big\}
\end{equation}
and $\|x_*^j\|=1$, $j=1,\ldots,s$. The normal cone in $(\mathcal{DMP})$ is understood in the sense of convex analysis
\begin{equation*}
N_C(x):=\big\{v\in\mathbb{R}^n\;\big|\;\langle v,y-x\rangle\le 0\;\textup{ for all }\;y\in C\big\}\;\textup{ if }\;x\in C,\,\textup{ and }\;N_C(x):=\emptyset\;\textup{ if }\;x\notin C.
\end{equation*}
It follows from the above definition that in ($\mathcal{DMP}$) we automatically have the {\em pointwise state constraints}
\begin{equation*}
 x(t)\in C\;\mbox{ for all }\;t\in[0,T].   
\end{equation*}
The constrained controlled sweeping system in ($\mathcal{DMP}$) is clearly a particular case of the time-delayed system ($\mathcal{P}$), while now we consider its optimization.

Throughout this and next sections, the following {\em standing assumptions} are imposed:

{\bf(H1)} The control set $U$ is a nonempty and compact subset of $\mathbb{R}^d$.

{\bf(H2)} The perturbation mapping $g:[0,T]\times\mathbb{R}^n\times\mathbb{R}^n\times\mathbb{R}^d\rightarrow\mathbb{R}^n$ is continuous in $(t, x,y,u)$ being also Lipschitz continuous with respect to $(t,x,y)$ uniformly on $U$ whenever $(x,y)$ belongs to a bounded subset of the product space $\mathbb{R}^n\times\mathbb{R}^n$.

{\bf(H3)} There exists a number $\beta>0$ ensuring the linear growth condition 
$$
\|g(t, x,y,u)\|\le \beta(1+\|x\|+\|y\|)\;\textup{ for all }\;(x,y,u)\in \mathbb{R}^n\times\mathbb{R}^n\times U.
$$

{\bf(H4)} The delay function $\delta:[0,T] \to [0,\infty)$ is Lipschitz continuous, nonincreasing and such that $0<\Delta:=\max\{\delta(t)\;|\;t\ge 0\}<\infty$.

{\bf(H5)} The history function $\phi:[-\Delta,0]\rightarrow\mathbb{R}^n$ is Lipschitz continuous.

It follows from Theorem~\ref{existence} that for each measurable control $u(\cdot)$, there exists a unique solution $x(\cdot )\in W^{1,2}([0, T],\mathbb{R}^n)$ to the  delay-differential inclusion in ($\mathcal{DMP}$). By {\em a feasible process} for ($\mathcal{DMP}$), we understand a pair $(x(\cdot),u(\cdot ))$ such that $u(\cdot)$ is measurable, $x(\cdot) \in W^{1,2}([0,T],\mathbb{R}^n)$, and all the constraints on controls and trajectories in ($\mathcal{DMP}$) are satisfied.

We are interested in developing a constructive scheme that allows us to approximate any feasible process for ($\mathcal{DMP}$) by feasible solutions to a perturbed time-discrete version of ($\mathcal{DMP}$). To furnish this, recall first the following result for an undelayed system taken from  \cite[Proposition~3.3]{B}.

\begin{lema}\label{Brez} Let $A:\mathbb{R}^n\longmapsto\mathbb{R}^n$ be a set-valued maximal monotone operator, let $f\in BV([0,T],\mathbb{R}^n)$, and let $x\in C([0,T],\mathbb{R}^n)$ be a solution to the differential inclusion
\begin{equation}\label{slow}
\dot x(t)\in-A(x(t))+f(t)\;\mbox{ a.e. }\;t\in [0,T].
\end{equation}
Then the following properties are equivalent:

{\bf(i)} We have $x(0)\in D(A)$ with $D(A)$ standing for the domain of $A$, i.e., the collection of points where the operator values are nonempty. 

{\bf(ii)} $x(\cdot)$ is Lipschitz continuous on $[0,T]$.\par
\noindent
In this case, $x(t)\in D(A)$ for all $t\in [0,T]$, and if we set 
$$
f(t^+):=\lim_{h\downarrow 0}\,\frac{1}{h}\int_{t}^{t+h} f(r)\,dr,\;t\in [0,T),
$$
then $x(\cdot)$ is right differentiable at every $t\in [0,T)$ with its right derivative $\dot x^+(\cdot)$ satisfying 
$$
\dot x^+(t)=-\textup{proj}_{A(x(t))}(f(t^+)) + f(t^+).
$$
If we assume in addition that $f\in W^{1,1}([0,T],\mathbb{R}^n)$, then $\dot x(\cdot)$ is right continuous at each $t\in [0,T)$. Furthermore, $\dot x(\cdot)$ is continuous at $t\in (0,T)$ if and only if $\dot x(\cdot)$ is continuous at $t$, and so $x(\cdot)$ is differentiable at this point. When moreover $A$ is single-valued and $f(\cdot)$ is continuous, then $x(\cdot)$ is weakly differentiable on $(0,T)$ and $\dot x(\cdot)$ is weakly continuous on $(0,T)$.
\end{lema}
\noindent
To formulate the next lemma, recall that a {\em representative} of a given measurable function on $[t_0,T]$ is a function (with certain prescribed properties) that agrees with the given one for a.e. $t\in[t_0,T]$.

\begin{lema}\label{auxiliar} Let $\phi(\cdot)$ be 
a history function, and let $(\overline{x}(\cdot),\overline{u}(\cdot))$ be an associated feasible pair for \textup{($\mathcal{DMP}$)} such that $\overline{u}(\cdot)$ is of bounded variation \textup{(}BV\textup{)} and admits a right continuous representative on $[0,T]$, which we keep denoting by $\overline{u}(\cdot)$. In addition to \textup{(H1)}, \textup{(H2)}, \textup{(H4)}, and \textup{(H5)}, suppose that the mapping $g(t,x,y,u)$ is locally Lipschitzian around $(t,\overline{x}(t),\overline{x}(t-\delta(t)),\overline{u}(t))$ for all $t\in [0,T]$. Then:

{\bf(i)} $\overline{x}(\cdot)$ is Lipschitz continuous on $[0,T]$ and right differentiable on this interval.

{\bf(ii)} The pair $(\overline{x}(t),\overline{u}(t))$, where the right derivative $\dot{\overline{x}}(t)$ is taken from {\rm(i)} and  where $\overline{u}(t)$ is a right continuous representative of $\overline{u}(\cdot)$, satisfies all the relationships in \textup{($\mathcal{DMP}$)} for all $t\in[0,T]$,
\end{lema}
{\bf Proof.} Denote $A(x):=N_C(x)$ and $f(t):=g(t, \overline{x}(t),\overline{x}(t-\delta(t)),\overline{u}(t))$. It is well known that $A(\cdot)$ is a maximal monotone multifunction. Moreover, the definition of $A(\cdot)$ implies that $D(A)=C$. Since $\overline{x}(\cdot)$ is absolutely continuous, $\overline{u}(\cdot)$ is of bounded variation, and the assumptions (H4) and 
(H5) are satisfied, it follows that the mapping $h(t):=(t, \overline{x}(t),\overline{x}(t-\delta(t)),\overline{u}(t))$ is of bounded variation on $[0,T]$, and hence its image is bounded in $[0,T]\times\mathbb{R}^n\times \mathbb{R}^n\times\mathbb{R}^d$. If $\Omega$ denotes the closure of this image, then $g$ is globally Lipschitz continuous on $\Omega$. Therefore, the composition $f(\cdot)=(g\circ h)(\cdot)$ is of bounded variation on $[0,T]$. Recall that $(\overline{x}(\cdot),\overline{u}(\cdot))$ is a feasible pair for \textup{($\mathcal{DMP}$)}. Thus
\begin{align}
\dot{\overline{x}}(t)&\in -A(\overline{x}(t))+f(t)\;\textup{ a.e. }\;t\in[0,T],\nonumber\\
\overline{x}(s)&=\phi(s),\;s\in[-\Delta,0],\label{cons1}\\
\overline{u}(t)&\in U,\;t\in [0,T],\label{cons2}
\end{align}
with the fulfillment of the state constraint $\overline{x}(0)\in D(A)$. It follows from Lemma~\ref{Brez} that $\overline{x}(\cdot)$ is Lipschitz continuous on $[0,T]$ and right differentiable on $[0,T)$, which justifies assertion (i). Using further the right continuity of $\overline{u}(\cdot)$, we get that $f(t^+)=f(t)$ for all $t\in [0,T)$. Employing again Lemma~\ref{Brez} gives us
$$
\dot{\overline{x}}^{\,+}(t)=-\textup{proj}_{A(\overline{x}(t))}(f(t)) + f(t)\in  -N_C(\overline{x}(t)) + g\big(t, \overline{x}(t),\overline{x}(t-\delta(t)),\overline{u}(t)\big)
$$
for all $t\in [0,T)$, which being combined with (\ref{cons1}) and (\ref{cons2}) justifies (ii) and thus completes the proof. \qed\vspace*{0.05in}

Now we are ready to establish the aforementioned central theorem, which provides a {\em well-posed} construction of {\em discrete approximations} of the time-delayed constrained sweeping process in ($\mathcal{DMP}$) and verifies the possibility to {\em strongly approximate} in the $W^{1,2}\times L^2$-norm topology of {\em any feasible solution} $(\overline{x}(\cdot),\overline{u}(\cdot))$ to ($\mathcal{DMP}$) by a sequence of feasible solutions to discrete-time sweeping processes that are properly extended to the continuous-time interval $[0,T]$. 

For any fixed $m\in\mathbb{N}$, consider the {\em uniform partition/mesh} on $[t_0,T]$  defined by
$$
\Delta_m:=\big\{0=t_0<t_1<\ldots < t_m^{2^m}=T\big\}\;\textup{ with }\;h_m:=t_m^{i+1} - t_m^i=\frac{T}{2^m}.
$$
\begin{theor}\label{convergence}
Let $\phi(\cdot)$ be a history function, and let $(\overline{x}(\cdot),\overline{u}(\cdot))$ be a feasible process for problem \textup{($\mathcal{DMP}$)} with a control representative $\overline{u}(\cdot)$ of bounded variation and right continuous on $[0,T)$. Imposing the standing assumptions above, we claim that for each $m\in\mathbb{N}$ there exists state-control pairs $(x_m(t),u_m(t))$ and perturbation terms $r_m(t)\ge 0$ and $\rho_m(t)\in \boldsymbol{\mathbb{B}}$ as $0\le t \le T$, for which the following hold:

{\bf(i)} The sequence of control mappings $u_m:[0,T]\rightarrow U$, which are constant on each interval $I_m^i:=[t_m^{i-1},t_m^i]$, converges to $\overline{u}(\cdot)$ strongly in $L^2([0,T],\mathbb{R}^d)$ and pointwise on $[0,T]$.

{\bf(ii)} The sequence of continuous state mappings $x_m:[0,T]\rightarrow\mathbb{R}^n$, which are affine on each interval $I_m^i$, converges strongly in $W^{1,2}([0,T],\mathbb{R}^n)$ to $\overline{x}(\cdot)$ and satisfy the discrete-time state constraints 
\begin{equation*}
x_m(t_m^i) =\overline{x}(t_m^i)\in C\;\;\textup{ for each }\;i=1,\ldots,2^m\;\textup{ with }\; x_m(0)=\phi(0)=x_0.
\end{equation*}

{\bf(iii)} For all $t\in(t_{m}^{i-1},t_m^i)$ and $i=1,\ldots,2^m$, we have 
\begin{align}
\dot x_m(t)&\in -N_C\big(x_m(t_m^i)\big)+g\big(t_m^i,x_m(t_m^i),x_m(t_m^i-\delta(t_m^i)),u_m(t)\big)+r_m(t)\rho(t),\label{discincl}\\
x_m(s)&=\phi(s),\;s\in[-\Delta,0],\label{discpast}\\
u_m(t)&\in U,\;t\in [0,T],\label{disccont}
\end{align}
where the mappings $r_m:[0,T]\rightarrow [0,\infty)$ and $\rho_m:[0,T]\rightarrow \boldsymbol{\mathbb{B}}$ are constant on each interval $I_m^i$ with 
\begin{equation*}
r_m(\cdot)\rightarrow 0\;\;\textup{ in }\;L^2(0,T)\;\;\textup{ as }\;m\rightarrow\infty. 
\end{equation*}
Moreover, all $x_m(\cdot)$ are Lipschitz continuous on $[0,T]$ with the same Lipschitz constant as $\ox(\cdot)$.
\end{theor} 
{\bf Proof.} For $m\in\mathbb{N}$ and $i=0,\ldots,2^m - 1$, define
\begin{align*}
x_m(t):&=\overline{x}(t_m^i) + (t-t_m^i)\frac{\overline{x}(t_m^{i+1})-\overline{x}(t_m^i)}{h_m},\;t\in [t_m^{i},t_m^{i+1}),\\
x_m(s):&=\phi(s),\;s\in [-\Delta,0],\\
u_m(t):&=\overline{u}(t_m^{i+1}),\;t\in [t_m^{i},t_m^{i+1}).
\end{align*}
Denote by $\omega_m(\cdot)$ the right derivative of $x_m(\cdot)$ on $[0,T)$ for which we have the representation
\begin{equation}\label{w}
\omega_m(t)=\omega_m^i:=\frac{\overline{x}(t_m^{i+1})-\overline{x}(t^i_m)}{h_m}\;\;\textup{ whenever }t\in [t_m^i, t^{i+1}_m),\;i=0,\ldots,2^m - 1.
\end{equation}

It follows from the right continuity of $\overline{u}(t)$ that $u_m(t)\rightarrow \overline{u}(t)$ and hence $\|u_m(t)-\overline{u}(t)\|^2\rightarrow 0$ as $m\rightarrow \infty$ for all $t\in [0,T)$. Moreover, the boundedness of $\overline{u}(\cdot)$ and the construction of $u_m(\cdot)$ ensure that the sequence of $\|u_m(\cdot)-\overline{u}(\cdot)\|^2$ is bounded by some constant $M>0$. The Lebesgue dominated convergence theorem tells us that $u_m(\cdot)\rightarrow \overline{u}(\cdot)$ strongly in $L^2(0,T)$ as $m\to\infty$, which justifies (i).

To proceed further, let $\overline{t}$ be a mesh point in the $m$-th partition of $[0,T]$. Since $h_m=T/2^m$ for each $m\in\mathbb{N}$, it follows that $\overline{t}$ remains a mesh point in the $m'$-th partition for any $m'\ge m$. Denoting  by $i_m(\overline{t})$ the index $i$ with $\overline{t}=i\frac{T}{2^m}$ and employing 
Lemma~\ref{auxiliar}(i) tell us that the right derivative $\dot{\overline{x}}(\overline{t})$ exists and
\begin{equation}\label{limder}
\lim_{m\to\infty}\,\omega_m^{i_m(\overline{t})}=\lim_{m\to\infty}\frac{\overline{x}\big(t_m^{i_m(\overline{t})+1}\big)-\overline{x}\big(t^{i_m(\overline{t})}_m\big)}{h_m}=\lim_{m\to\infty}\frac{\overline{x}(\overline{t}+h_m)-\overline{x}(\overline{t})}{h_m}=\dot{\overline{x}}(\overline{t}).
\end{equation}
We are going to show that the sequence $\omega_m(\cdot)$ converges to $\dot{\overline{x}}(\cdot)$ in $L^2((0,T),\mathbb{R}^n)$. Recall from 
Lemma~\ref{auxiliar} that $\overline{x}(\cdot)$ is Lipschitz, and hence its derivative $\dot{\overline{x}}(t)$ exists for a.e. $[0,T]$ by the seminal Rademacher theorem. Let us verify that $\omega_m(t)\to \dot{\overline{x}}(t)$ for a.e. $t\in[0,T]$, where $\omega_m(\cdot)$ are taken from \eqref{w}. To furnish this, denote by $\tau_m(t)$ the unique mesh point $t_m^i$ such that $t\in [t_m^i,t_m^{i+1})$, and then set $\alpha_m(t):=t_m^{i+1}-t$ and $\beta_m(t):=t-t_m^i$. Observe that $\tau_m(t)+h_m=t+(t_m^{i+1}-t)=t+\alpha_m(t)$ and that $0<\alpha_m(t)\le h_m$. Similarly we get $\tau_m(t)=t-\beta_m(t)$ and $0<\beta_m(t)\le h_m$. If $t$ belongs to the full measure subset of $[0,T]$ where $\dot{\overline{x}}(t)$ exists (hence both $\dot{\overline{x}}(t^+)$ and $\dot{\overline{x}}(t^-)$ exist and agree with each other), then 
\begin{align*}
\omega_m(t)&=\frac{1}{h_m}\big(\overline{x}(t_{i+1}^m)-\overline{x}(t^i_m)\big)\\
&=\frac{1}{h_m}\bigg(\frac{\overline{x}(\tau_m(t)+h_m)-\overline{x}(t)}{\tau_m(t)+h_m-t}(\tau_m(t)+h_m-t) + \frac{\overline{x}(t)-\overline{x}(\tau_m(t))}{t-\tau_m(t)}(t-\tau_m(t))\bigg)\\
&=\frac{1}{h_m}\bigg(\frac{\overline{x}(t+\alpha_m(t))-\overline{x}(t)}{\alpha_m(t)}(\tau_m(t)+h_m-t) + \frac{\overline{x}(t)-\overline{x}(t-\beta_m(t))}{\beta_m(t)}(t-\tau_m(t))\bigg)\\
&=\frac{1}{h_m}\Big((\dot{\overline{x}}(t^+)+o(1))(\tau_m(t)+h_m-t) + (\dot{\overline{x}}(t^-)+o(1))(t-\tau_m(t))\Big)\\
&=\frac{1}{h_m}(\dot{\overline{x}}(t)+o(1))h_m=\dot{\overline{x}}(t)+o(1)\rightarrow \dot{\overline{x}}(t)\;\mbox{ as }\;m\to\infty,
\end{align*}
and therefore $\|\omega_m(t)-\dot{\overline{x}}(t)\|^2\rightarrow 0$ a.e. $t\in[0,T]$. Using again dominated convergence theorem yields
\begin{equation}\label{derconv}
 \lim_{m\to \infty}\,\|\omega_m(\cdot)-\dot{\overline{x}}(\cdot)\|_{L^2((0,T);\mathbb{R}^n)}=0,   
\end{equation}
which verifies the claim. If $L>0$ is a Lipschitz constant of $\overline{x}(\cdot)$ on $[0,T]$, then for every $m\in\mathbb{N}$ and each $t\in [0,T]$ we get the estimates
\begin{align*}
\|x_m(t) - \overline{x}(t)\|&\le  \|\overline{x}(t_m^i) - \overline{x}(t)\| + \frac{|t-t_m^i|}{h_m}\|\overline{x}(t_m^{i+1}) - \overline{x}(t_m^i)\|\\
&\le L|t_m^i - t| + L|t_m^{i+1} - t|\\
&\le 2Lh_m\rightarrow 0,\;\textup{ as }\;m\to\infty.
\end{align*}
The latter ensures the uniform convergence $\|x_m(\cdot) - \overline{x}(\cdot)\|_{\infty}^2\rightarrow 0$ on $[0,T]$, and hence the strong convergence of $x_m(\cdot)$ to $\overline{x}(\cdot)$ in $W^{1,2}([0,T],\mathbb{R}^n)$ in view of (\ref{derconv}), which thus justifies (ii).

It remains to verify (iii). Fixing $m\in\mathbb{N}$ and employing 
Lemma~\ref{auxiliar}(ii) together with (\ref{limder}) at every mesh point $\overline{t}$ tell us that
$$
\lim_{m\to\infty}\,\omega_m^{i_m(\overline{t})}=\dot{\overline{x}}(\overline{t})\in -N_C\big(\overline{x}(\overline{t})\big)+g\big(\overline{t},\overline{x}(\overline{t}),\overline{x}(\overline{t}-\delta(\overline{t})),\overline{u}(\overline{t})\big).
$$
Let $\overline{\zeta}\in N_C(\overline{x}(\overline{t}))$ be such that $\dot{\overline{x}}(\overline{t})+\overline{\zeta}-g(\overline{t},\overline{x}(\overline{t}),\overline{x}(\overline{t}-\delta(\overline{t})),
\overline{u}(\overline{t}))=0$. Recall that from certain $m\in\mathbb{N}$ onwards we have $x_m(\overline{t})=\overline{x}(\overline{t})$. The Lipschitz continuity of $\delta(\cdot)$, the uniform convergence of $x_m(\cdot)$ to $\overline{x}(\cdot)$, both on $[0,T]$, and the fact that $x_m(s)=\overline{x}(s)=\phi(s)$ for all $s\in [-\Delta,0]$, imply that $\lim_{m\to\infty} x_m(\overline{t}-\delta(\overline{t}))=\overline{x}(\overline{t}-\delta(\overline{t}))$. Additionally, the right continuity of $\overline{u}(\cdot)$ ensures that $\lim_{m\to\infty} u_m(\overline{t})=\lim_{m\to\infty} \overline{u}(\overline{t}+h_m)= \overline{u}(\overline{t})$. The latter observations and the continuity of $g$ yield
$$
\lim_{m\to\infty}\,\big(\omega_m^{i_m(\overline{t})}+\overline{\zeta}-g(\overline{t},x_m(\overline{t}),x_m(\overline{t}-\delta(\overline{t})),
u_m(\overline{t}))\big)=\dot{\overline{x}}(\overline{t})+\overline{\zeta}-g(\overline{t},\overline{x}(\overline{t}),\overline{x}(\overline{t}-
\delta(\overline{t})),\overline{u}(\overline{t}))=0.
$$
Thus there exists a sequence $\{r_m(\overline{t})\}$, with  $r_m(\overline{t})\downarrow 0$ as $m\to\infty$ and such that 
$$
\omega_m^{i_m(\overline{t})}\in-N_C(x_m(\overline{t}))+g(\overline{t},x_m(\overline{t}),x_m(\overline{t}-\delta(\overline{t})),u_m(\overline{t}))
+r_m(\overline{t})\mathbb{B}.
$$
The previous inclusion allows us to choose a vector $\rho_m(\overline{t})\in\mathbb{B}$ satisfying
$$
\omega_m^{i_m(\overline{t})}=-N_C\big(x_m(\overline{t}))+g(\overline{t},x_m(\overline{t}),x_m(\overline{t}-\delta(\overline{t})),u_m(\overline{t})\big)
+r_m(\overline{t})\rho_m(\overline{t}).
$$
Finally, we extend $r_m(\cdot)$ and $\rho_m(\cdot)$ to the 
entire interval $[0,T]$ by setting $r_m(t):=r_m(t_m^i)$ and $\rho_m(t):=\rho_m(t_m^i)$ when $t\in [t_m^i,t_m^{i+1})$. Recall that $\dot x_m(t)=\omega_m^{i_m(\overline{t})}$ whenever $t\in [t_m^{i_m(\overline{t})},t_m^{i_m(\overline{t})+1})$. Therefore,
$$
\dot x_m(t)\in -N_C\big(x_m(t_m^i))+g(t_m^i,x_m(t_m^i),x_m(t_m^i -\delta(t_m^i)),u_m(t_m^i)\big)+r_m(t)\rho_m(t),\;\;t\in (t_m^i,t_m^{i+1}),
$$
which is inclusion (\ref{discincl}). Parts (\ref{discpast}) and (\ref{disccont}) result clear from the definitions of $x_m(\cdot)$ on $[-\Delta,0]$ and $u_m(\cdot)$, respectively. The remaining part of the theorem is a direct consequence of the constructions.\qed

\section{Discrete Approximations of Sweeping Local Minimizers}

The goal of this section is to verify the strong $W^{1,2}\times L^2$-approximation of {\em local minimizers} in the 
time-delayed sweeping optimal control problem $(\mathcal{DMP})$ by optimal solutions to its discrete-time counterparts. First we introduce the following appropriate notion of local minimizers for $(\mathcal{DMP})$.

\begin{defi}\label{localmin} A feasible pair $(\overline{x}(\cdot ),\overline{u}(\cdot))$ for $(\mathcal{DMP})$ is a $W^{1,2}\times L^2$-local minimizer in this problem if there exists $\epsilon > 0$ such that $J[\overline{x},\overline{u}]\leq J[x,u]$ for all
feasible pairs $(x(\cdot),u(\cdot))$ satisfying
$$
\int _0^T \Bigl(
\| \dot x(t)-\dot{\overline{x}}(t)\|^2+\|u(t)-\overline{u}(t)\|^2\Bigr)dt < \epsilon.
$$
\end{defi}
\noindent
This notion, which is an extension of the one introduced in \cite{cmn18} for undelayed system, obviously covers a larger territory than its {\em strong} version, where the $W^{1,2}$-norm for $x(\cdot)$ is replaced by the ${\cal C}$-norm.

To proceed with the constructions of discrete approximation of local minimizers for the time-delayed sweeping optimal control problem $(\mathcal{DMP})$, we need to modify the notion from Definition~\ref{localmin} by involving its {\em relaxation} of the Bogolyubov-Young type, which relates to a certain convexification with respect to velocities and has been well understood in the calculus of variations and optimal control; see, e.g., \cite{m06}. 

\begin{defi}\label{relaxed} A feasible pair $(\overline{x}(\cdot ),\overline{u}(\cdot))$ for $(\mathcal{DMP})$ is a relaxed $W^{1,2}\times L^2$-local minimizer in this problem if there exists $\epsilon > 0$ such that
\begin{equation*}
\ph\big(\ox(T)\big)\le\ph\big(x(T)\big)\;\textrm{ whenever }\;\int_0^{T}\left(\n\dot{x}(t)-\dot{\ox}(t)\en^2+\n u(t)-\ou(t)\en^2\right)dt<\e,
\end{equation*}
where $u(\cdot)$ is a measurable control satisfying the convexified constraint $u(t)\in\co U$  for a.e. 
$t\in[0,T]$, and  where $x(\cdot)$ is a trajectory of the convexified inclusion
\begin{equation}\label{cog}
\dot x(t)\in-N_C\big(x(t)\big)+\textup{co}\;g\big(t,x(t),x(t-\delta(t)),U\big),
\end{equation}
which can be uniformly approximated in $W^{1,2}([0,T];\R^n)$-norm topology by feasible trajectories to $(\mathcal{DMP})$ generated by piecewise constant controls $u^k (\cdot)$ on $[0,T]$ in such a way that their convex combinations
strongly converges to $u(\cdot)$ in the norm topology $L^2([0,T];\R^d)$.
\end{defi}\vspace*{-0.1in}
\noindent
There is obviously no difference between the local minimizers from Definitions~\ref{localmin} and \eqref{relaxed} under the convexity assumptions imposed on the sets $U$ and $g(\cdot,U)$. On the other hand, the beauty and power of relaxation for continuous-time differential systems, as first observed by Bogolyubov and Young in the classical calculus of variations, is that the above {\em local relaxation stability} holds automatically without convexity assumptions, mainly in the case of strong local minimizers. We refer the reader to \cite{cmn18} and the bibliographies therein for more discussions in the case of undelayed systems.

Using the approximation results of Theorem~\ref{convergence}, we can construct a suitable 
sequence of discrete-time optimal control problems whose optimal solutions strongly converge to a designated local minimizer of the original time-delayed sweeping control problem $\mathcal{DMP}$. This process not only ensures a robust connection between the discrete and continuous formulations but also enables us to approximate the local minimizer in question by optimal solutions for the discrete-time counterparts. To furnish this, fix a relaxed $W^{1,2}\times L^2$-local minimizer $(\overline{x}(\cdot),\overline{u}(\cdot))$ of the optimal control problem ($\mathcal{DMP}$), take $\epsilon$ from 
Definition~\ref{localmin}, and for any $m\in\mathbb N$ define the multifunction
\begin{equation*}
 F_m(t,x,y,u):=N\big(x(t);C)-g(t,x(t),y(t),u(t)\big)-r_m(t)\rho_m(t),
\end{equation*}
where $y(t):= x(t-\delta(t))$ on $t\in [0,T]$.
We construct the following family of discrete-time problems ($\mathcal{P}_m$), where $r_m(\cdot )$ and $\rho_m(\cdot )$ are taken from the formulation of 
Theorem~\ref{convergence}:
\begin{equation}\label{d_a_p}
\begin{array}{ll}
&\mbox{minimize}\;\; J_m[x_m,u_m]:=\varphi(x_m^{2^m})\\
&+\disp\sum_{i=0}^{2^m-1}\int_{t_m^i}^{t_m^{i+1}}\left(\n \frac{x_m^{i+1}-x_m^i}{h_m}-\dot{\bar{x}}(t)\en^2+\n \frac{y_m^{i+1}-y_m^i}{h_m}-\dot{\bar{x}}(t-\delta(t))\en^2+\n u_m^i-\bar{u}(t)\en^2\right)dt
\end{array}
\end{equation}
 over the collection of the discrete functions 
 $$
 (x_m,u_m):=(x_m^0,x_m^1,\ldots, x_m^{2^m},u_m^0,u_m^1\ldots,u_m^{2^m-1})\;\mbox{ and }\;y^i_m:= x_m(t^i_m-\delta(t^i_m))
 $$
for $i=0,\ldots,2^m-1$ satisfying the geometric and functional constraints given by
\begin{equation}\label{re_1}
x_m^{i+1}- x_m^i\in -h_mF_m(t_m^i,x_m^i,y_m^i, u_m^i)\;
\textrm{ for }\;i=0,\ldots,2^m-1,
\end{equation}
\begin{equation}
\label{re_6}
\la x^j_\ast, x_m^{2^m}\ra \leq c_j\textrm{ for all }\;j=1,\ldots,s,
\end{equation}
\begin{equation*}
x_m^0:=\phi(0)=x_0\in C,
\end{equation*}
\begin{equation*}\label{re_2}
\left\|\frac{x_m^{i+1}-x_m^i}{h_m}\right\|\le L\;\text{ for all } i=0,\ldots , 2^m-1,
\end{equation*}
\begin{equation}\label{re_5}
\| x_m^i-\ox(t_m^i)\|\le\e/2\; \textrm{ for } i=0,\ldots,2^m-1,
\end{equation}
\begin{equation}\label{re_7}
\| y_m^i-\ox(t_m^i-\delta(t_m^i))\|\le\e/2\; \textrm{ for }\;i=0,\ldots,2^m-1,
\end{equation}
\begin{equation*}\label{re_4}
u_m^i\in U\;\textrm{ for }\;i=0,\ldots,2^m-1,
\end{equation*}
\begin{equation}\label{re_3}
\sum_{i=0}^{2^m-1}\int_{t_m^i}^{t_m^{i+1}}\left(\n\frac{x_m^{i+1}-x_m^i}{h_m}-\dot{\ox}(t)\en^2+\n \frac{y_m^{i+1}-y_m^i}{h_m}-\dot{\bar{x}}(t-\delta(t))\en^2+\n u_m^i-\ou(t)\en^2\right)dt\le \e/2.
\end{equation}

To proceed further, we have to ensure that the discrete problems ($\mathcal{P}_m$) admit optimal solutions.
 
\begin{prop}
In addition to the standing assumptions, suppose that the cost $\varphi$ is lower semicontinuous on $\mathbb{R}^n$. Then for each problem ($\mathcal{P}_m$) there exists an optimal solution provided that $m\in\mathbb N$ is sufficiently large.   
\end{prop}
{\bf Proof.} Since any relaxed $W^{1,2}\times L^2$-local minimizer $(\overline{x}(\cdot),\overline{u}(\cdot))$ is a feasible solution to ($\mathcal{DMP}$) by Definition~\ref{relaxed}, we can apply to the designated minimizer the results of Theorem~\ref{convergence}. This gives us a sequence of triples 
$(x_m(\cdot),y_m(\cdot),u_m(\cdot))$ with $y_m(t):=x_m(t-\delta(t))$ that satisfy the assertions (i)--(iii) of that theorem and ensures, in particular, the fulfillment of the constraints in \eqref{re_5}, \eqref{re_7}, and \eqref{re_3} for all large $m\in\mathbb N$. Therefore, the set of feasible solutions to ($\mathcal{P}_m$) is a nonempty subset of $\mathbb{R}^{2^m +1}\times\mathbb{R}^{2^m +1}\times\mathbb{R}^{2^m}$ for large $m$. Since the perturbation $g(\cdot,\cdot,\cdot,\cdot)$, delay $\delta(\cdot)$, and history $\phi(\cdot)$ are all continuous functions while $C$ and $U$ are closed and compact sets, respectively, we get that the limit of any convergent sequence of feasible solutions to ($\mathcal{P}_m$) is a feasible solution itself, and hence the set of feasible solutions is closed for each fixed $m\in\mathbb{N}$. Let us now verify that the latter set is bounded.

To furnish this, we deduce from (\ref{re_3}) the relationships
\begin{equation}\label{expand}
\begin{array}{ll}
&\disp\n\frac{x_m^{i+1}-x_m^i}{h_m}\en^2 h_m-2\int_{t_m^i}^{t_m^{i+1}}\disp\Big\langle \frac{x_m^{i+1}-x_m^i}{h_m},\disp\dot{\ox}(t)\Big\rangle\,dt +\disp\int_{t_m^i}^{t_m^{i+1}}\disp\n\dot{\ox}(t)\en^2 dt\\
&=\disp\int_{t_m^i}^{t_m^{i+1}}\n\frac{x_m^{i+1}-x_m^i}\disp{h_m}-\dot{\ox}(t)\en^2 dt<\e/2.
\end{array}
\end{equation}
Considering further the nonnegative numbers
$$
\beta:=\int_{0}^{T}\n\dot \ox(t)\en dt,\;\;\; p_i=\n x_m^{i+1}-x_m^i\en,\;\;i=0,1,\ldots,2^{m}-1,
$$
and using the above conditions (\ref{expand}) tell us that
$$
p_i^2 < 2\beta p_i + \frac{\e}{2}h_m,\;\;i=0,1,\ldots,2^{m}-1,
$$
which immediately implies the inequality
$$
(p_i-\beta)^2 < \beta^2+\frac{\e}{2}h_m,\;\;i=0,1,\ldots,2^{m}-1.
$$
Denote $\gamma:=(\beta^2 + (\e/2)h_m)^{1/2}+\beta$, $r_i:=\|x_m^i\|$ and deduce from the above inequality that
$$
r_{i+1}< r_i + \gamma,\;\;i=0,1,\ldots,2^m-1. 
$$
Applying \cite[Exercise~4.1.8]{CLSW} on each index subset $0\le j\le i-1$ with $\delta_j=0$ and $\Delta_j=\gamma$ yields $\|x_m^i\|\le 2^m\gamma$ for all $i=0,1,\ldots,2^m$.  Similarly we can derive a bound for $u_m^i$. Then the continuity of the functions $\delta(\cdot)$, $\oy(\cdot)$, and $\phi(\cdot)$ combined with (\ref{re_5}) yields the claimed boundedness. Taking into account the lower semicontinuity of the cost function in \eqref{d_a_p}, we finally deduce the existence of optimal solutions from the Weierstrass existence theorem in finite-dimensional spaces. $\h$

Now we are ready to establish the main result of this section on the $W^{1,2}\times L^2$-strong convergence of optimal solutions in the discrete problems $(\mathcal{P}_m)$ to the designated local minimizer of the original time-delayed sweeping control problem ($\mathcal{DMP}$).

\begin{theor} Let $\phi(\cdot)$ be a history function, and let $(\ox(\cdot ),\ou(\cdot ))$ be a relaxed $W^{1,2} \times L^2$-local minimizer for problem $(\mathcal{DMP})$. Consider any sequence $\{(\ox_m(\cdot),\ou_m(\cdot))\}$ of optimal solutions to problems $(\mathcal{P}_m)$, $m\in\mathbb N$ and extend $\ox_m(\cdot)$ piecewise linearly while $\ou_m(\cdot)$ piecewise constantly to the interval $[0,T]$ without changing the labels. Suppose in addition to the standing assumptions that $\varphi$ is continuous around $\ox(T)$. Then we get the strong convergence
$$
\bigl(\ox_m(\cdot),\ou_m(\cdot)\bigr)
\rightarrow
\big((\ox(\cdot),\ou(\cdot)\big)\;\mbox{ as }\;m\to\infty
$$
in the norm topology of $W^{1,2}([0,T];\mathbb{R}^n)\times L^2([0,T];\mathbb{R}^d)$.
\end{theor}
{\bf Proof}. We start with establishing the inequality
\begin{equation}\label{J}
\disp\limsup_{m\to\infty}J_m[\ox_m,\ou_m]\le J[\ox,\ou]
\end{equation}
for any sequence of optimal solutions to $(P_m)$. To verify \eqref{J}, suppose while arguing by for contradiction, that there exists a sequence of natural numbers $m\to\infty$ and a positive number $\zeta$ satisfying
\begin{equation}\label{cont_J}
J[\ox,\ou]=\varphi(\ox(T))<J_m[\ox_m,\ou_m]-\zeta\;\text{ for all }\;m\in\N.
\end{equation}
Construct by Theorem~\ref{convergence} an approximating sequence $\bigl(
x_m(\cdot),u_m(\cdot)\bigr)$ and deduce from the continuity of the original cost function $\varphi$ around $x(T)$ that
$$
\varphi(x_m(T))\to \varphi(\ox(T))\;\mbox{ as }\; m\to \infty.
$$
It follows from the approximating cost expression \eqref{d_a_p} along the $(x_m(\cdot), u_m(\cdot))$, due to the $W^{1,2}\times L^2$-strong convergence $(x_m(\cdot), u_m(\cdot
))\to (\ox(\cdot),\ou(\cdot
))$ by Theorem~\ref{convergence} and the continuity of $\delta(\cdot)$, that 
\begin{equation*}
\begin{aligned}
& \sum_{i=0}^{2^m-1}\int_{t^i_m}^{t^{i+1}_m}\Bigg(\n\frac{x^{i+1}_m-x^i_m}{h_m}-\dot{\ox}(t)\en^2+\n \frac{y_m^{i+1}-y_m^i}{h_m}-\dot{\bar{x}}(t-\delta(t))\en^2+\n u^i_m-\ou(t)\en^2\Bigg)dt \\
&=\int^{T}_{0}\Big(\n\dot{x}_m(t)-\dot{\ox}(t)\en^2+\n\dot{x}_m(t-\delta(t))-\dot{\ox}(t-\delta(t))\en^2+\n u^i_m-\ou(t)\en^2\Big)dt\to 0
\end{aligned}
\end{equation*}
as $m\to \infty$, where $x_m(s):=\phi(s)$ for $s\in[-\Delta, 0]$ taken from \eqref{discpast}. This ensures that
$$
J_m[x_m(\cdot),u_m(\cdot)]\rightarrow J[\ox,\ou]\;\;\mbox{as}\;\;m
\rightarrow\infty.
$$
The latter convergence combined with the feasibility of $(x_m(\cdot), u_m(\cdot))$ to problem $(\mathcal{P}_m)$ contradicts \eqref{cont_J} when $m$ is large enough, and thus \eqref{J} is justified.

Our next step is to show that
\begin{equation}\label{e:6.10}
\underset{m\to\infty}{\mathrm{lim}}\Big[\gamma_m:=\int_0^T\Big(\n\dot{\ox}_m(t)-\dot{\ox}(t)\en^2+\n\dot{y}_m(t)-\dot{\ox}(t-\delta(t))\en^2+\n\ou_m(t)-\ou(t)\en^2\Big)dt\Big]=0
\end{equation}
which clearly yields the convergence of $(\ox_m(\cdot), \ou_m(\cdot))$ to $(\ox(\cdot),\ou(\cdot))$ in the norm topology of $W^{1,2}([0,T],\mathbb{R}^n)\times L^2([0,T], \mathbb{R}^d)$. 
Arguing by contradiction, suppose that the limit of $\gamma_m$ in \eqref{e:6.10}, along a subsequence without relabeling, equals some number $\gg>0$. It follows from \eqref{re_3} that the sequence of extended optimal solutions $\{(\dot\ox_m(\cdot),\ou_m(\cdot))\}$ to $(P_m)$ with $\oy_m(t):=\ox_m(t-\delta(t))$ satisfy
\begin{equation*}
\sum_{i=0}^{2^m-1}\int_{t_m^i}^{t_m^{i+1}}\left(\n\frac{\ox_m^{i+1}-\ox_m^i}{h_m}-\dot{\ox}(t)\en^2+\n \frac{\oy_m^{i+1}-\oy_m^i}{h_m}-\dot{\bar{x}}(t-\delta(t))\en^2+\n \ou_m^i-\ou(t)\en^2\right)dt\le\e/2,
\end{equation*}
}
which ensures that $\{(\dot\ox_m(\cdot),\ou_m(\cdot))\}$ is bounded in the space $L^2([0,T];\R^n)\times L^2([0,T];\R^d)$. Due to the weak compactness of the unit ball in reflexive space, there exist functions $(\Tilde v(\cdot),\Tilde u(\cdot))\in L^2([0,T];\R^n)\times L^2([0,T];\R^d)$ such that a subsequence of $\{(\dot\ox_m(\cdot),\ou_m(\cdot))\}$  (without relabeling) weakly converges to $(\Tilde v(\cdot),\Tilde u(\cdot))$ in the product space $W^{1,2}\times L^2$. Define further the function
\begin{equation*}
\tilde{x}(t):= \begin{cases}
x_0+\disp\int_{t_0}^t\tilde v(\tau)d\tau\;\textrm{ for all }\;t\in[t_0,T],\\
\phi(t) \;\textrm{ for all }\;t\in[-\Delta, t_0],
\end{cases}
\end{equation*}
which gives us $\dot{\tilde x}(t)=\tilde v(t)$ for a.e. $t\in[0,T]$ and implies that the functions $\dot\ox_m(\cdot)$ converge weakly to $\tilde v(\cdot)= \dot\tx(\cdot)$ in $L^2([0,T];\R^{n})$. Furthermore, Mazur's weak closure theorem ensures the existence of a sequence of convex combinations of $\big(\dot\ox_m(\cdot),\ou_m(\cdot)\big)$ that converges strongly to $\big(\dot{ \tilde x}(\cdot),\tilde u(\cdot)\big)$ in $W^{1,2}([0,T];\R^n)\times L^2([0,T];\R^d)$ and thus almost everywhere on 
$[0,T]$ along their subsequences (without relabeling). The obtained pointwise convergence of convex combinations allows us to conclude that $\tilde u(t)\in \text{co}U$ for a.e. $t\in[0,T]$ and $(\tx(\cdot), \tu(\cdot))$ satisfies the differential inclusion \eqref{cog}. 

Consider now the integral functional 
$$
\mathcal{I}[\vartheta]:=\int_0^T\|\vartheta(t)-\dot\ox(t)\|^2dt,
$$
which is lower semicontinuous in the $L^2$-weak topology by the convexity of the integrand in $y$. Hence 
\begin{equation*}
\begin{array}{ll}
\disp\mathcal{I}[\dot\tx, \disp\tu]=\int^T_0\left\|\left(\dot\tx(t),\dot\tx(t-\delta(t)), \disp\tu(t)\right)-\left(\dot \ox(t),\dot\ox(t-\disp\delta(t)),\ou(t)\right)\right\|^2dt\\
\label{I}\leq\disp\liminf_{m\to\infty}\sum_{i=0}^{2^m-1}\disp\int_{t^i_m}^{t^{i+1}_m}
\n\disp\bigg(\dfrac{\ox^{i+1}_m-\ox^i_m}{h_m},\frac{\oy_m^{i+1}-\oy_m^i}{h_m},\ou^i_m\bigg)-(\dot\ox(t),\dot\ox(t-\delta(t)), \ou(t))
\en^2dt
\end{array}
\end{equation*}
by the construction of $\tx(\cdot)$. Taking into account the continuity of $\delta(\cdot)$, passing to the limit in \eqref{re_5}, \eqref{re_7}, and \eqref{re_3} as $m\to \infty$ and then using \eqref{I} yields
$$
\|\tx(t)-\ox(t)\|\leq \e/2,\quad\|\tx(t-\delta(t))-\ox(t-\delta(t))\|\leq \e/2,
$$
$$
\int_0^T\|(\dot\tx(t),\dot\tx(t-\delta(t)),\tu(t))-(\dot\ox(t),\dot\ox(t-\delta(t))\ou(t))\|^2dt \leq \e/2, 
$$
which ensures that the pair $(\tx(\cdot), \tu(\cdot))$ belongs to the prescribed $W^{1,2}\times L^2$-neighborhood of the designated local minimizer $(\ox(\cdot), \ou(\cdot ))$ and satisfies the relationships
$$
J[\tx,\tu]+\gg/2= \varphi(\tx(T))+\gg/2\leq\liminf_{m\to\infty}\left[\varphi(\ox_m^{2^m})+\gg/2 \right] = \liminf_{m\to\infty}J_m[\ox_m, \ou_m]\leq J[\ox,\ou]
$$
by the continuity of $\varphi$ around $\ox(T)$,  \eqref{J}, and \eqref{d_a_p}. This contradicts the fact that $(\ox(\cdot), \ou(\cdot)$ is a relaxed $W^{1,2}\times L^2$-local minimizer for $(\mathcal{P}_m)$ and thus completes the proof of this theorem. $\h$

\section{Conclusions}

The paper establishes well-posedness and stability with respect to discrete approximations for a new class of optimal control problems governed by time-delayed discontinuous sweeping differential inclusions. The obtained strong approximation of feasible solutions by finite-dimensional systems with discrete time is of both theoretical and numerical interest providing potentials for further algorithmic developments. On the other hand, the derived strong convergence of optimal solutions for discrete-time problems to a local minimizer of the time-delayed controlled sweeping process makes a bridge between delayed continuous-time and discrete-time systems. In our future research, we plan to utilize the obtained results, married to appropriate tools of variational analysis and generalized differentiation, in deriving necessary optimality conditions foir discrete approximation problems and then, by passing to the limit over the diminishing step of discretization, for local minimizers of delayed continuous-time sweeping control systems.

\end{document}